\documentclass[11pt]{article}
\usepackage{epsfig}
\def\be{\begin{equation}}
\def\ee{\end{equation}}
\def\bea{\begin{eqnarray}}
\def\eea{\end{eqnarray}}
\def\bes{\begin{eqnarray*}}
\def\ees{\end{eqnarray*}}

\def\nn{\nonumber}
\def\<{\langle}
\def\>{\rangle}
\def\lb{\label}
\def\bs{\setminus}

\def\R{{\bf R}}
\def\C{{\bf C}}
\def\Z{{\bf Z}}
\def\N{{\bf N}}
\def\U{{\bf U}}
\def\Q{{\bf Q}}
\def\T{{\bf T}}
\def\CP{{\bf CP}}

\def\aa{{\alpha}}
\def\bb{{\beta}}
\def\ga{{\gamma}}

\def\th{{\theta}}
\def\om{{\omega}}
\def\Om{{\Omega}}
\def\ep{{\epsilon}}
\def\lm{{\lambda}}
\def\Lm{{\Lambda}}
\def\dl{{\delta}}
\def\sg{{\sigma}}

\def\Sg{{\Sigma}}
\def\vf{{\varphi}}

\def\vep{{\varepsilon}}

\def\A{{\cal A}}
\def\H{{\cal H}}
\def\T{{\cal T}}

\def\P{{\cal P}}

\def\Nn{{\cal N}}

\def\im{{\rm im}}

\def\Sp{{\rm Sp}}

\def\dm{{\rm \diamond}}

\def\ol#1{\overline{#1}}  %overline in math mode

\def\hb{\vrule height0.18cm width0.14cm $\,$}

\def\ol#1{\overline{#1}}  %overline in math mode

\title{Resonance identity, stability and multiplicity of \\
closed characteristics on compact convex hypersurfaces}

\author{\thanks{This paper has been accepted by
{\it Duke Math. J. }}    Wei Wang$^{1}$,\thanks{ Partially supported
by NNSF and RFDP of MOE of China. E-mail:
alexanderweiwang@yahoo.com.cn }\qquad Xijun
Hu$^{3}$,\thanks{Partially supported by NNSF of China (No.
10526038). E-mail: xjhu@amss.ac.cn }\qquad Yiming
Long$^{1,2}$,\thanks{Partially supported by the 973 Program of MOST,
Yangzi River Professorship, NNSF, MCME, RFDP, LPMC of MOE of
China, and Nankai University. E-mail: longym@nankai.edu.cn}  \\ \\
$^{1}$ Chern Institute of Mathematics\\
$^{2}$ Key Lab of Pure Mathematics and Combinatorics of Ministry of Education\\
Nankai University\\
Tianjin 300071, The People's Republic of China\\
$^{3}$ Institute of Mathematics, Academy of Mathematics and Systems Science\\
Chinese Academy of Sciences, Beijing 100080, The People's Republic of China\\}

\date{}

\begin{document}

\maketitle

\begin{abstract}
{\it There is a long standing conjecture in Hamiltonian analysis
which claims that there exist at least $n$ geometrically distinct
closed characteristics on every compact convex hypersurface in
$\R^{2n}$ with $n\ge 2$. Besides many partial results, this
conjecture has been only completely solved for $n=2$. In this
paper, we give a confirmed answer to this conjecture for $n=3$. In
order to prove this result, we establish first a new resonance
identity for closed characteristics on every compact convex
hypersurface $\Sg$ in $\R^{2n}$ when the number of geometrically
distinct closed characteristics on $\Sg$ is finite. Then using
this identity and earlier techniques of the index iteration
theory, we prove the mentioned multiplicity result for $\R^6$. If
there are exactly two geometrically distinct closed
characteristics on a compact convex hypersuface in $\R^4$, we
prove that both of them must be irrationally elliptic. }
\end{abstract}



{\bf Running head}: Closed characteristics on convex hypersurfaces

\renewcommand{\theequation}{\thesection.\arabic{equation}}
\renewcommand{\thefigure}{\thesection.\arabic{figure}}

\setcounter{equation}{0}%\setcounter{figure}{0}
\section{Introduction and main results}%{Section 1}

In this paper, let $\Sigma$ be a fixed $C^3$ compact convex hypersurface
in $\R^{2n}$, i.e., $\Sigma$ is the boundary of a compact and strictly
convex region $U$ in $\R^{2n}$. We denote the set of all such hypersurfaces
by $\H(2n)$. Without loss of generality, we suppose $U$ contains the origin.
We consider closed characteristics $(\tau, y)$ on $\Sigma$, which are
solutions of the following problem
\be
\left\{\matrix{\dot{y}=JN_\Sigma(y), \cr
               y(\tau)=y(0), \cr }\right. \lb{1.1}\ee
where $J=\left(\matrix{0 &-I_n\cr
                I_n  & 0\cr}\right)$,
$I_n$ is the identity matrix in $\R^n$, $\tau>0$, $N_\Sigma(y)$ is
the outward normal vector of $\Sigma$ at $y$ normalized by the
condition $N_\Sigma(y)\cdot y=1$. Here $a\cdot b$ denotes the
standard inner product of $a, b\in\R^{2n}$. A closed characteristic
$(\tau, y)$ is {\it prime}, if $\tau$ is the minimal period of $y$.
Two closed characteristics $(\tau, y)$ and $(\sigma, z)$ are {\it
geometrically distinct},  if $y(\R)\not= z(\R)$. We denote by
$\T(\Sg)$ the set of all geometrically distinct closed
characteristics on $\Sg$. A closed characteristic $(\tau,y)$ is
{\it non-degenerate}, if $1$ is a Floquet multiplier of $y$ of precisely
algebraic multiplicity $2$, and is {\it elliptic}, if all the Floquet
multipliers of $y$ are on ${\bf U}=\{z\in\C\,|\,|z|=1\}$, i.e., the unit
circle in the complex plane.

There is a long standing conjecture on the number of closed characteristics
on compact convex hypersurfaces in $\R^{2n}$:
\be \,^{\#}\T(\Sg)\ge n, \qquad \forall \; \Sg\in\H(2n). \lb{1.2}\ee

Since the pioneering works \cite{Rab1} of P. Rabinowitz and \cite{Wei1}
of A. Weinstein in 1978 on the existence of at least one closed
characteristic on every hypersurface in $\H(2n)$, the existence of multiple
closed characteristics on $\Sg\in\H(2n)$ has been deeply studied
by many mathematicians. When $n\ge 2$, besides many results under pinching
conditions, in 1987-1988 I. Ekeland-L. Lassoued, I. Ekeland-H. Hofer, and
A, Szulkin (cf. \cite{EkL1}, \cite{EkH1}, \cite{Szu1}) proved
$$ \,^{\#}\T(\Sg)\ge 2, \qquad \forall\,\Sg\in\H(2n). $$
In \cite{LoZ1} of 2002, Y. Long and C. Zhu further proved
$$ \;^{\#}\T(\Sg)\ge [\frac{n}{2}]+1, \qquad \forall\, \Sg\in \H(2n), $$
where we denote by $[a]\equiv\max\{k\in\Z\,|\,k\le a\}$. Note that this
estimate yields still only at least $2$ closed characteristics when $n=3$.
We refer readers to the survey paper \cite{Lon5} and the recent \cite{Lon6} of
Y. Long for earlier works and references on this conjecture. Our following
main result in this paper gives a confirmed answer to the conjecture
(\ref{1.2}) for $n=3$.

{\bf Theorem 1.1.} {\it There holds $\;^{\#}\T(\Sg)\ge 3$ for every
$\Sg\in\H(6)$. }

One of the main ingredients of our proof of this theorem is a new
resonance identity on closed characteristics. In \cite{Eke1} of
1984, I. Ekeland discovered that there must exist a resonance
condition relating the closed characteristics on $\Sg\in\H(2n)$
provided $\,^{\#}\T(\Sg)<+\infty$. However, he did not state
explicitly what the resonance condition is. Then in \cite{Vit1} of
1989, C. Viterbo clarified such a resonance condition by
establishing a mean index identity for closed characteristics on
compact star-shaped hypersurfaces in $\R^{2n}$ provided all closed
characteristics on $\Sg$ together with their iterations are
non-degenerate (cf. p.234 of \cite{Eke3}). Note that in \cite{Rad1}
of 1989 and \cite{Rad2} of 1992, a similar identity for closed
geodesics on compact Finsler manifolds was established by H.-B.
Rademacher. Motivated by these results, in the current paper we
establish the following mean index identity for closed
characteristics on every $\Sg\in\H(2n)$ when
$\,^{\#}\T(\Sg)<+\infty$. This yields hopefully an explicit version
of what I. Ekeland discovered.

{\bf Theorem 1.2.} {\it Suppose $\Sigma\in \H(2n)$ satisfies
$\,^{\#}\T(\Sg)<+\infty$. Denote all the geometrically distinct closed
characteristics by $\{(\tau_j,\; y_j)\}_{1\le j\le k}$. Then the following
identity holds
\be  \sum_{1\le j\le k}\frac{\hat{\chi}(y_j)}{\hat{i}(y_j)}=\frac{1}{2},
\lb{1.3}\ee
where $\hat{i}(y_j)\in\R$ is the mean index of $y_j$ given by Definition 3.14,
$\hat{\chi}(y_j)\in\Q$ is the average Euler characteristic given by
Definition 3.15 and Remark 3.16 below. Specially by (\ref{3.56}) below we
have
\be \hat\chi(y) = \frac{1}{K(y)}\sum_{1\le m\le K(y)\atop 0\le l\le 2n-2}
  (-1)^{i(y^{m})+l}k_l(y^{m}),  \lb{1.4}\ee
$K(y)\in \N$ is the minimal period of critical modules of
iterations of $y$ defined in Proposition 3.13, $i(y^{m})$ is
the Morse index of a corresponding dual-action functional at the
$m$-th iteration $y^m$ of $y$ (cf. Definition 3.3 and Proposition
3.5 below), $k_l(y^{m})$ is the critical type numbers of $y^m$
given by Definition 3.11 below. }

{\bf Remark 1.3.} Note that $1/2$ in the right hand side of (\ref{1.3})
comes from the average Euler characteristic for equivariant homology on the
space of loops in $\R^{2n}$. In fact, we have
$$ \frac{1}{2}
=\lim_{N\rightarrow\infty}\frac{1}{N}\sum_{q\le N}\dim H_q(\CP^\infty;\Q)
=\lim_{N\rightarrow\infty}\frac{1}{N}\sum_{q\le N}\dim H_q(S^{\infty}\times_{S^1}\{0\};\Q). $$
Since the space of loops in $\R^{2n}$ is $S^1$-equivariantly homotopic to its origin
$\{0\}$, the last term in the above expression represents the average Euler
characteristic of the $S^1$-equivariant homology of the space of loops in $\R^{2n}$.
For details, we refer to Section 5 below.

When all the closed characteristics on $\Sigma\in\H(2n)$
together with their iterations are nondegenerate, by Remark 3.16 below our
identity (\ref{1.3}) coincides with the identity (1.3) of Theorem 1.2 of
\cite{Vit1} as well as $(153)$ in p.234 of \cite{Eke3}. Thus for
$\Sg\in\H(2n)$ our Theorem 1.2 generalizes C. Viterbo's result in \cite{Vit1}
to degenerate closed characteristics.

Note that in \cite{HWZ1} of 1998, H. Hofer-K. Wysocki-E. Zehnder proved
that $\,^{\#}\T(\Sg)=2$ or $\infty$ holds for every $\Sg\in\H(4)$.
In \cite{Lon3} of 2000, Y. Long proved further that $\Sg\in\H(4)$ and
$\,^{\#}\T(\Sg)=2$ imply that both of the closed characteristics must be
elliptic, i.e., each of them possesses four Floquet multipliers with two
$1$s and the other two locate on the unit circle too. Now as a by-product
of our Theorem 1.2 we obtain a stronger result:

{\bf Theorem 1.4.} {\it Let $\Sg\in\H(4)$ satisfy
$\,^{\#}\T(\Sg)=2$. Then both of the closed characteristics must
be irrationally elliptic, i.e., each of them possesses four
Floquet multipliers with two $1$s and the other two located on the
unit circle with rotation angles being irrational multiples of
$\pi$. }

Because of above mentioned results and other indications, we suspect that the
following conjectures hold:

{\bf Conjecture 1.5.} {\it For every integer $n\ge 2$, there holds }
$$ \{\,^{\#}\T(\Sg)\,|\,\Sg\in {\cal H}(2n)\} = \{n\}\cup\{+\infty\}. $$

It seems that for $n\ge 3$ there is no effective methods so far which can be
used to prove that $\,^{\#}\T(\Sg)>n$ implies $\,^{\#}\T(\Sg)=\infty$.

Recall that a closed characteristic is irrationally elliptic, if it is elliptic
and the linearized Poincar\'e map is suitably homotopic to the $\dm$-product of
one $\left(\matrix{1 & 1 \cr
                   0 & 1 \cr}\right)$ and $n-1$ rotation $2\times 2$ matrices
with rotation angles being irrational multiples of $\pi$. Note that based upon our
studies on the stabilities of closed characteristics on $\Sg\in {\cal H}(2n)$,
and closed geodesics on Finsler spheres, we tend to believe that the following
may hold.

{\bf Conjecture 1.6.} {\it All the geometrically distinct closed characteristics
on $\Sg$ are irrationally elliptic for $\Sg\in {\cal H}(2n)$ with $n\ge 2$
whenever $\,^{\#}\T(\Sg)<\infty$. }

The rest of this paper is arranged as follows.

$\<1\>$ Motivated by the works \cite{Kli1} and \cite{Kli2} of W. Klingenberg,
\cite{GrM2} of D. Gromoll and W. Meyer, \cite{Eke1} and \cite{Eke3} of
I. Ekeland and \cite{Vit1} of C. Viterbo, for every $\Sg\in\H(2n)$ with
$\,^{\#}\T(\Sg)<+\infty$, we shall construct a functional $\Psi_a$ for large
$a>0$ on the space of loops in $\R^{2n}$ and establish a Morse theory of this
functional $\Psi_a$ to study closed characteristics on $\Sg$.

As usual we use the Clarke-Ekeland dual action principle and a modification
of the Ekeland index theory. Because in general such a dual action functional
is not $C^2$, motivated by the studies on closed geodesics and convex
Hamiltonian systems, we follow \cite{Eke1} to introduce a finite dimensional
approximation to the space of loops in $\R^{2n}$ to get the enough smoothness.
For the dual action principle, the
origin becomes an accumulation point of its critical values. To estimate the
contribution of the critical point at the origin to the Morse Series, we
construct a special family of Hamiltonian functions which have better properties
at the origin and infinity, and are homogenous in the middle. Such a construction
allows us to give a precise understanding of the behavior of the dual action
functional near the origin.

$\<2\>$ In Section 2, fixing a hypersurface $\Sg\in\H(2n)$ with
$\,^{\#}\T(\Sg)<+\infty$, we construct a family $\A$ of Hamiltonian functions
by Proposition 2.4 using auxiliary functions satisfying conditions (i)-(iv) of
Proposition 2.2. Using such Hamiltonian functions, we construct a functional
$\Psi_a$ on the space of loops in $\R^{2n}$ for every $a>0$ whose critical
points are precisely all the closed characteristics on $\Sg$ with periods
less than $a$ and that the origin of the loop space is the only constant
critical point of $\Psi_a$.

$\<3\>$ In Section 3, we prove that for every fixed closed characteristic
$(\tau,y)$ on $\Sg$, the critical modules of all the functionals $\Psi_a$
produced by $H_a\in\A$ at its critical point corresponding to $(\tau,y)$ are
isomorphic to each other whenever $a>\tau$. Therefore we can further require
the Hamiltonian function in $\A$ to be homogeneous near such  critical points
so that the critical modules  are periodic functions of the
dimension. This homogeneity of the Hamiltonian function is realized by the
condition (v) of Proposition 2.2.

$\<4\>$ Using the properties of the Hamiltonian functions in $\A$, in Section
4, the property of the dual action functional near the origin is understood
precisely and we show that the origin has in fact no homological contribution
to the lower order terms in the Morse series.

$\<5\>$ Using the homological information obtained in the Sections 2-4, in
Section 5, we compute all the local critical modules of the dual action
functional $\Psi_a$ and use such information to set up a Morse theory for
all the closed characteristics on $\Sg\in\H(2n)$. Together with the global
homological information on the loop space we establish the claimed mean index
identity (\ref{1.3}) and prove Theorem 1.2.

$\<6\>$ Using Theorem 1.2 together with the techniques developed in the index
iteration theory we give proofs of Theorems 1.1 and 1.4 in Section 6.

Here we give a brief sketch for the proof of Theorem 1.1. Assuming that Theorem
1.1 does not hold, i.e., $\,^{\#}\T(\Sg)\le 2$, by \cite{EkH1} or \cite{LoZ1}
we should have $\,^{\#}\T(\Sg)=2$. By our Theorem 1.2 the two prime closed
characteristics $(\tau_1,y_1)$ and $(\tau_2,y_2)$ must satisfy the following
identity:
\be \frac{\hat{\chi}(y_1)}{\hat{i}(y_1)}
         + \frac{\hat{\chi}(y_2)}{\hat{i}(y_2)}=\frac{1}{2}.  \lb{1.5}\ee
Here it is well known that $\hat{i}(y_j)>2$ for $j=1$ and $2$ always holds
(cf. Theorem 1.7.7 of \cite{Eke3} or Lemma 15.3.2 of \cite{Lon4}).
By \cite{LoZ1} (cf. Theorem 15.5.2 of \cite{Lon4}), at least one of $\hat{i}(y_1)$
and $\hat{i}(y_2)$ is irrational, say $\hat{i}(y_1)\in \R\bs\Q$ without loss of
generality.

Then if $\hat{\chi}(y_1)\not= 0$, we obtain $\hat{i}(y_2)$ is irrational too by
(\ref{1.5}). Thus by a careful derivation using the index iteration formulae of
Y. Long in \cite{Lon3} and estimates obtained by Y. Long and C. Zhu in \cite{LoZ1},
there should exist more than two closed characteristics which yields a
contradiction.

If $\hat{\chi}(y_1)= 0$, by the index iteration formulae of Y. Long in
\cite{Lon3}, one can prove that the orbit $y_2$ must satisfy
\be \frac{\hat{\chi}(y_2)}{\hat{i}(y_2)}\le \frac{1}{4}. \lb{1.6}\ee
Then together with (\ref{1.5}), it yields a contradiction too and proves the
theorem.

In this paper, let $\N$, $\N_0$, $\Z$, $\Q$, $\R$, and $\R^+$ denote
the sets of natural integers, non-negative integers, integers, rational
numbers, real numbers, and positive real numbers respectively.
Denote by $a\cdot b$ and $|a|$ the standard inner product and norm in
$\R^{2n}$. Denote by $\langle\cdot,\cdot\rangle$ and $\|\cdot\|$
the standard $L^2$ inner product and $L^2$ norm. For an $S^1$-space $X$, we
denote by $X_{S^1}$ the homotopy quotient of $X$ by $S^1$, i.e.,
$X_{S^1}=S^\infty\times_{S^1}X$, where $S^\infty$ is the unit sphere in an infinite
dimensional {\it complex} Hilbert space. In this paper we use only $\Q$
coefficients for all homological modules. By $t\to a^+$, we mean $t>a$ and
$t\to a$.

\setcounter{equation}{0}%\setcounter{figure}{0}
\section{A variational structure for closed characteristics}%{Section 2}

In the rest of this paper, we fix a $\Sg\in\H(2n)$ and assume the following condition
on $\T(\Sg)$:

\noindent (F) {\bf There exist only finitely many geometrically distinct
closed characteristics \\$\quad \{(\tau_j, y_j)\}_{1\le j\le k}$ on $\Sigma$. }

In this section, we transform the problem (\ref{1.1}) into a fixed period
problem of a Hamiltonian system and then study its variational structure.
We introduce the following set:

{\bf Definition 2.1.} {\it Under the assumption (F), the set of
periods on $\Sigma$ is defined by }
$$ per(\Sigma)=\{m\tau_j\;|\; m\in\N,\; 1\le j\le k\}.$$

Clearly $per(\Sigma)$ is a discrete subset of $\R^+$. Motivated by
Definition 2.1 of \cite{Eke1} and Lemma 2.2 of \cite{Vit1}, we construct
the following auxiliary function to further define Hamiltonian functions.

{\bf Proposition 2.2.} {\it  For any sufficiently small $\vartheta\in(0,1)$,
there exists a function $\varphi\equiv \varphi_{\vartheta}\in C^\infty(\R, \;\R^+)$
depending on $\vartheta$ which has $0$ as its unique critical point in $[0, +\infty)$
such that the following hold

(i) $\varphi(0)=0=\varphi^\prime(0)$ and
$ \varphi^{\prime\prime}(0)=1=\lim_{t\rightarrow 0^+}\frac{\varphi^\prime(t)}{t}$.

(ii) $\varphi(t)$ is a polynomial of degree $2$ in a neighborhood of $+\infty$.

(iii) $\frac{d}{dt}\left(\frac{\varphi^\prime(t)}{t}\right)<0$ for $t>0$,
and $\lim_{t\rightarrow +\infty}\frac{\varphi^\prime(t)}{t}<\vartheta$,
i.e., $\frac{\varphi^\prime(t)}{t}$ is strictly decreasing for $t> 0$.

(iv) $\min(\frac{\varphi^\prime(t)}{t}, \varphi^{\prime\prime}(t))\ge \sigma$
for all $t\in \R^+$ and some $\sigma>0$. Consequently, $\varphi$ is strictly
convex on $[0,\,+\infty)$.

(v) In particular, we can choose $\aa\in (1,2)$ sufficiently close to $2$
and $c\in (0,1)$ such that $\vf(t)=ct^{\aa}$ whenever
$\frac{\vf'(t)}{t}\in [\vartheta, 1-\vartheta]$ and $t>0$. }

{\bf Proof. } We construct a $\varphi$ satisfying (i)-(v).

Define $\varphi_1(t)=\left(\frac{\alpha^2-7\alpha+12}{2}\right)t^2
+(-\alpha^2+6\alpha-8)t^3
+\left(\frac{\alpha^2-5\alpha+6}{2}\right)t^4 $ for  $t\in (-\infty,
1]$ and $\varphi_1(t)=t^\alpha$ for $t\in [1, +\infty)$. Then
$\varphi_1\in C^2$ and $\frac{\varphi_1^\prime(t)}{t}$ is strictly
decreasing for $t>0$ and $\alpha\in(1,\, 2)$. Note that
$\lim_{t\rightarrow 0^+}\frac{\varphi_1^\prime(t)}{t}
=\varphi_1^{\prime\prime}(0)=\alpha^2-7\alpha+12>2$ and
$\varphi_1^\prime(1)=\alpha<2$. Hence
$\varphi_2(t)\equiv\frac{\varphi_1(t)}{\varphi_1^{\prime\prime}(0)}$
satisfies (i) and (iii). Next we further modify $\varphi_2$ in a
neighborhood of $+\infty$.

Since $\frac{\varphi_2^\prime(t)}{t}$ tends to $0$ when $t$
goes to $+\infty$, we obtain a $T=T_\vartheta>1$ sufficiently large
such that
\be \frac{\varphi_2^\prime(t)}{t}=c\alpha t^{\alpha-2}
\in \left(0, \frac{\vartheta}{2\aa-1}\right),\quad \forall t\ge T, \lb{2.1}\ee
where $c=c_{\aa}=(\aa^2-7\aa + 12)^{-1}$. Now we define
\be \varphi_3(t)=\left\{\matrix{ \varphi_2(t),\quad &&t\in (-\infty, T]; \cr
cT^\alpha+c\alpha T^{\alpha-1}(t-T)+
\frac{c\alpha(\alpha-1) T^{\alpha-2}}{2}(t-T)^2,
\quad &&t\in[T, +\infty). \cr }\right.
\lb{2.2}\ee
Then $\varphi_3\in C^\infty (\R\setminus\{1, T\}, \R^+)\cap C^2(\R, \R^+)$.
We can approximate $\varphi_3$ by a smooth function $\varphi$
(cf. Theorem 2.5 of \cite{Hir1}) such that $\varphi=\varphi_3$ holds
outside a small neighborhood of $\{1, T\}$, and
$\|\varphi-\varphi_3\|_{C^2}$ is small enough. Then it is easy to
see that $\varphi$ satisfies (i)-(iv) of the proposition.

Note that $f(t)\equiv\frac{\varphi_2^\prime(t)}{t}$ is a strictly decreasing
function with $f(1)=\frac{\alpha}{\alpha^2-7\alpha+12}<1$. Since $f(1)$ tends
to $1$ as $\alpha$ goes to $2$, if $\alpha$ is chosen sufficiently close
to $2$, then $f(t)>1-\vartheta$ holds for all $t\in [0,1]$. Together with (\ref{2.1})
it is easy to verify that this $\vf$ satisfies (v). \hfill\hb

{\bf Remark 2.3.} $1^{\circ}$ Note that in the above proof for (v), we can choose
$\aa\in (1,2)$ sufficiently close to $2$, $c\in (0,1)$, and $T$
sufficiently large such that $\vf(t)=ct^{\aa}$ if and only if $1+\dl\le t\le T-\dl$
for some $\dl>0$. Note that the property (v) above is used only in the second part
of the proof of Proposition 3.5 below to show that our index and nullity given by
Definition 3.3 below coincide with those defined in \cite{Eke1}-\cite{Eke3}, and
in our study in the Subsection 3.2 to obtain the periodic property of critical modules
at critical points. In the other parts of this paper we use functions $\vf$ which
satisfy the properties (i)-(iv).

$2^{\circ}$ Note that in the following, we only need the definition of $\varphi$ on
$[0,\,+\infty)$. In Proposition 2.4 below, the parameter $\vartheta$ given by
Proposition 2.2 depends on the parameter $a$, i.e., given an $a>\hat{\tau}$ as in
Proposition 2.4, we choose first the parameter
$\vartheta\in (0,\frac{1}{a}\min\{\hat{\tau},\hat{\sg}\})$. Then we can choose the
parameter $\alpha\in(1,\,2)$ depending on $a$
and let $\varphi$ to be homogeneous of degree $\alpha$ and then modify it near $0$ and
$+\infty$ such that (i)-(iv) in Proposition 2.2 hold. Here we do not require $\varphi$
to satisfy (v) in Proposition 2.2. We denote
such choices of $\vartheta$, $\aa$ and $\varphi$ by $\vartheta_a$, $\aa_a$ and
$\varphi_a$ respectively to indicate their dependence on $a$. In such a way, we
can obtain a connected family of $\varphi_a$ continuously depending on $a$. Each
$\varphi_a$ in this family satisfies properties (i)-(iv) of Proposition 2.2.
Moreover, the first and second derivatives of $\varphi_a(t)$ with
respect to $t$ are also continuous in the parameter $a$. Note that under these
choices, the coefficients of the polynomials in the proof of (ii) of Proposition 2.2
are continuous in $a$. Here that $\varphi_a$s form a connected family in $a$ is
crucial in our study below, for example in the proofs of Proposition 3.2, Lemma 3.4
and Proposition 3.5.

Let $j: \R^{2n}\rightarrow\R$ be the gauge function of $\Sigma$, i.e.,
$j(\lambda x)=\lambda$ for $x\in\Sigma$ and $\lambda\ge0$, then
$j\in C^3(\R^{2n}\setminus\{0\}, \R)\cap C^0(\R^{2n}, \R)$
and $\Sigma=j^{-1}(1)$. Denote by $\hat{\tau}=\inf\{s\,|\, s\in per(\Sigma)\}$
and $\hat{\sg}=\min\{|y|^2\,|\, y\in\Sigma\}$.

{\bf Proposition 2.4.} {\it Let $a>\hat{\tau}$,
$\vartheta_a\in\left(0,  \frac{1}{a}\min\{\hat{\tau}, \hat{\sg}\}\right)$
and $\varphi_a$ be a $C^\infty$ function associated to $\vartheta_a$ satisfying
(i)-(iv) of Proposition 2.2. Define the Hamiltonian function $H_a(x)=a\varphi_a(j(x))$
and consider the fixed period problem
\bea
\left\{\matrix{\dot{x}(t)=JH_a^\prime(x(t)) \cr
     x(0)=x(1)         \cr }\right. \lb{2.3}\eea
Then the following hold:

(i) $H_a\in C^3(\R^{2n}\setminus\{0\}, \R)\cap C^1(\R^{2n}, \R)$
and there exist $R, r>0$ such that
$$r|\xi|^2\le H^{\prime\prime}_a(x)\xi\cdot\xi\le R|\xi|^2,\quad
\forall x\in \R^{2n}\setminus\{0\},\;\xi\in \R^{2n}.$$

(ii) There exist $\epsilon_1, \epsilon_2\in \left(0, \frac{1}{2}\right)$
and $C\in \R$, such that
$$\frac{\epsilon_1 |x|^2}{2}-C\le H_a(x)\le\frac{\epsilon_2 |x|^2}{2}+C,
\quad \forall x\in\R^{2n}.$$

(iii) Solutions of (\ref{2.3}) are $x\equiv0$ and $x=\rho y(\tau t)$ with
$\frac{\varphi_a^\prime(\rho)}{\rho}=\frac{\tau}{a}$,
where $(\tau, y)$ is a solution of (\ref{1.1}).
In particular, nonzero solutions of (\ref{2.3})
are in one to one correspondence with solutions of (\ref{1.1})
with period $\tau<a$.

(iv) There exists $r_0>0$ independent of $a$ and there
exists $\mu_a>0$ depending on $a$ such that}
$$ H_a^{\prime\prime}(x)\xi\cdot\xi\ge {2ar_0}|\xi|^2,
    \qquad{\rm for}\quad 0<|x|\le \mu_a,\; \xi\in\R^{2n}. $$

{\bf Proof.} Since $j(\lambda x)=\lambda j(x)$ for all
$x\in \R^{2n}\setminus \{0\}$ and $\lambda\in \R^+$, we have
$j^\prime(\lambda x)=j^\prime(x)$ and
$j^{\prime\prime}(\lambda x)=\lambda^{-1}j^{\prime\prime}(x)$.
Hence for $x=\lambda y$ with $y\in\Sigma$ and $\lambda\in \R^+$ we have
\bea H_a^{\prime\prime}(x)\xi\cdot\xi&=&H_a^{\prime\prime}(\lambda y)\xi\cdot\xi\nn\\
&=&a\varphi_a^{\prime\prime}(j(\lambda y))(j^\prime(\lambda y)\cdot\xi)^2
+a\varphi_a^{\prime}(j(\lambda y))j^{\prime\prime}(\lambda y)\xi\cdot\xi\nn\\
&=&a\varphi_a^{\prime\prime}(\lambda)(j^\prime(y)\cdot\xi)^2
+a\varphi_a^{\prime}(\lambda)\lambda^{-1}j^{\prime\prime}(y)\xi\cdot\xi\lb{2.4}\\
&\ge& a\sigma((j^\prime(y)\cdot\xi)^2+j^{\prime\prime}(y)\xi\cdot\xi),\lb{2.5}
\eea
where the last inequality follows from (iv) of Proposition 2.2.
Now fix  $y\in\Sigma$ and represent $\R^{2n}=\R y\oplus T_y\Sigma$.
We define a new norm in $\R^{2n}$ by
\be |z|_y^2\equiv\lambda^2+|z_2|^2,\quad \forall
z=\lambda y+z_2\in\R y\oplus T_y\Sigma. \lb{2.6}\ee
Since any two norms on $\R^{2n}$ are equivalent, we have
\be C_1(y)^{-1}|z|\le |z|_y\le C_1(y)|z|, \lb{2.7}\ee
for some constant $C_1(y)>0$ depending on $y$. Note that
$j^\prime(y)=N_\Sigma(y)$ by the fact that $N_\Sigma(y)\cdot y=1$
and $j^\prime(y)\cdot y=j(y)=1$ for every $y\in\Sigma$.
Since $j(\lambda y)=\lambda j(y)$, we have $j^\prime(y)\cdot y=j(y)$,
hence $j^{\prime\prime}(y)y=0$. For $\xi=\lambda y+\xi_2$, we have
\bea (j^\prime(y)\cdot\xi)^2+j^{\prime\prime}(y)\xi\cdot\xi
&=&(j^\prime(y)\cdot(\lambda y+\xi_2))^2
+j^{\prime\prime}(y)(\lambda y+\xi_2)\cdot(\lambda y+\xi_2)\nn\\
&=&(j^\prime(y)\cdot\lambda y)^2+j^{\prime\prime}(y)\xi_2\cdot\xi_2\nn\\
&\ge& \lambda^2+ C_2(y)|\xi_2|^2\nn\\
&\ge& C_3(y)|\xi|^2,\nn
\eea
for some positive constants $C_2(y)$ and $C_3(y)$ depending on $y$.
Here the first inequality holds since $\Sigma$ is strictly convex,
hence $j^{\prime\prime}(y)|_{T_y \Sigma}$ is positive definite.
The last inequality follows from (\ref{2.6}) and (\ref{2.7}).
By the compactness of $\Sigma$ and (\ref{2.7}) we have
$H^{\prime\prime}_a(x)\xi\cdot\xi\ge r|\xi|^2$ for some $r>0$.
The compactness of $\Sigma$ and (\ref{2.4})  yield
$H^{\prime\prime}_a(x)\xi\cdot\xi\le R|\xi|^2$ for some $R>0$.
This proves (i).

For (ii), it suffices to consider $|x|$ large. Hence suppose
$x=\lambda y$ for $y\in\Sigma$ and $\lambda>0$, then $|x|=\lambda |y|$,
so $\lambda$ is large. By (ii) and (iii) of Proposition 2.2, we have
$H_a(x)=a\varphi_a(\lambda)=D_0+D_1\lambda+D_2\lambda^2$
for some $0<2D_2<a\vartheta$. Hence
$H_a(x)=D_0+\frac{D_1}{|y|}|x|+\frac{D_2}{|y|^2}|x|^2$
with $\frac{D_2}{|y|^2}<\frac{a\vartheta}{2|y|^2}<\frac{1}{2}$ by
the definition of $\hat{\sg}$. This proves (ii).

Clearly $x\equiv 0$ is the unique constant solution of (\ref{2.3}).
Suppose $x(t)$ is a nonconstant solution of (\ref{2.3}), then
$H_a(x(t))=a\varphi_a(j(x(t)))=const$. Since $\varphi_a$ is strictly
increasing, we have  $j(x(t))=const$. Let $\rho=j(x(t))$ and
$y(t)=\rho^{-1}x\left(\frac{\rho}{a\varphi_a^\prime(\rho)}t\right)$.
Then $j(y)=\rho^{-1}j(x)=\rho^{-1}\rho=1$, hence $y(\R)\subset \Sigma$.
Moreover, we have $\dot y(t)=JN_\Sigma(y(t))$ by (\ref{2.3}).
Hence $\left(\frac{a\varphi_a^\prime (\rho)}{\rho},\; y\right)$
is a solution of (\ref{1.1}). By (i) and (iii) of Proposition 2.2,
we have $\frac{\varphi^\prime_a(\rho)}{\rho}<1$. Hence
$\tau\equiv\frac{a\varphi^\prime_a(\rho)}{\rho}<a$.
This together with $a\vartheta_a<\hat{\tau}$ proves one side of (iii).
The other side of (iii) can be proved similarly and thus is omitted.

(\ref{2.4}) together with the proof of (i) and Proposition 2.2 (i)
yield (iv). \hfill\hb

In the following, we will use the Clarke-Ekeland dual action principle.
As usual, the Fenchel transform of a function
$F: \R^{2n}\rightarrow\R$ is defined by
\be F^\ast (y)=\sup\{x\cdot y-F(x)\;|\; x\in \R^{2n}\}. \lb{2.8}\ee
Following Proposition 2.2.10 of \cite{Eke3}, Lemma 3.1 of \cite{Eke1}
and the fact that $F_1\le F_2\Leftrightarrow F_1^\ast\ge F_2^\ast$, we
have:

{\bf Proposition 2.5.} {\it Let $H_a$ be a function
defined in Proposition 2.4 and $G_a=H_a^\ast$ the Fenchel
transform of $H_a$. Then we have

(i) $G_a\in C^2(\R^{2n}\setminus\{0\}, \R)\cap C^1(\R^{2n}, \R)$
and $$G_a^\prime(y)=x\Leftrightarrow y=H_a^\prime(x)\Rightarrow
H_a^{\prime\prime}(x)G_a^{\prime\prime}(y)=1.$$

(ii) $G_a$ is strictly convex. Let $R$ and $r$ be the real
numbers given by (i) of Proposition 2.4. Then we have
$$R^{-1}|\xi|^2\le G^{\prime\prime}_a(y)\xi\cdot\xi\le  r^{-1}|\xi|^2,
\quad \forall y\in \R^{2n}\setminus\{0\},\;\xi\in \R^{2n}.$$

(iii) Let $\epsilon_1, \epsilon_2, C$ be the real numbers given by (ii)
of Proposition 2.4. Then we have
$$\frac{|x|^2}{2\epsilon_2}-C\le G_a(x)\le\frac{|x|^2}{2\epsilon_1}+C,
\quad \forall x\in\R^{2n}.$$

(iv) Let $r_0>0$ be the constant given by (iv) of Proposition 2.4.
Then there exists $\eta_a>0$ depending on $a$ such that the following holds
$$G_a^{\prime\prime}(y)\xi\cdot\xi\le \frac{1}{2ar_0}|\xi|^2,
\qquad{\rm for}\quad 0<|y|\le \eta_a,\; \xi\in\R^{2n}. $$

(v) In particular, let $H_a=a\varphi_a(j(x))$ with $\varphi_a$ satisfying
further (v) of Proposition 2.2. Then we have
$G_a(\mu j^\prime(z))=c_1\mu^\beta$ when $z\in \Sigma$
and $\mu j^\prime(z)\in \{H_a^\prime (x)\;|\; H_a(x)=acj(x)^\alpha\}$,
where $c$ is given by (v) of Proposition 2.2, $c_1>0$ is some constant and
$\alpha^{-1}+\beta^{-1}=1$ holds with $\aa=\aa_a$ and $\bb=\bb_a$ depending
on $a$.} \hfill\hb

Now we apply the dual action principle to problem (\ref{2.3}). Let
\be L_0^2(S^1, \R^{2n})=
\left\{u\in L^2([0, 1],\R^{2n})\left|\frac{}{}\right.\int_0^1u(t)dt=0\right\}.
\lb{2.9}\ee
Define a linear operator
$M: L_0^2(S^1, \R^{2n})\rightarrow L_0^2(S^1, \R^{2n})$ by
\be \frac{d}{dt}Mu(t)=u(t), \quad \int_0^1Mu(t)dt=0.\lb{2.10}\ee
The dual action functional on $L_0^2(S^1, \R^{2n})$ is defined by
\be \Psi_a(u)=\int_0^1\left(\frac{1}{2}Ju\cdot Mu+G_a(-Ju)\right)dt,
\lb{2.11}\ee
where $G_a$ is given by Proposition 2.5.

By (ii) of Proposition 2.5 and the proof of Proposition 3.3
on p.33 of \cite{Eke1}, we have

{\bf Proposition 2.6.} {\it The functional $\Psi_a$
is $C^{1, 1}$ on $L_0^2(S^1, \R^{2n})$. Suppose $x$
is a solution of (\ref{2.3}), then $u=\dot{x}$ is a critical
point of $\Psi_a$. Conversely, suppose $u$ is a critical point
of $\Psi_a$, then there exists a unique $\xi\in\R^{2n}$ such that
$Mu-\xi$ is a solution of (\ref{2.3}). In particular, solutions
of (\ref{2.3}) are in one to one correspondence with critical
points of $\Psi_a$. }\hfill\hb

{\bf Proposition 2.7.} {\it The functional $\Psi_a$ is bounded from
below on $L_0^2(S^1, \R^{2n})$.}

{\bf Proof.} For any $u\in L_0^2(S^1, \R^{2n})$, we represent
$u$ by its Fourier series
\bea u(t)=\sum_{k\neq 0} e^{k2\pi Jt}x_k, \quad x_k\in\R^{2n}.\lb{2.12}\eea
Then we have
\bea Mu(t)=-J\sum_{k\neq 0} \frac{1}{2\pi k}e^{k2\pi Jt}x_k.\lb{2.13}\eea
Hence
\bea \frac{1}{2}\langle Ju,\; Mu\rangle=
-\frac{1}{2}\sum_{k\neq 0} \frac{1}{2\pi k}|x_k|^2
\ge -\frac{1}{4\pi}\|u\|^2.\lb{2.14}
\eea
By (\ref{2.11}), we have
\bea \Psi_a(u)&=&\int_0^1\left(\frac{1}{2}Ju\cdot Mu+G_a(-Ju)\right)dt\nn\\
&\ge&\frac{1}{2}\langle Ju,\; Mu\rangle
+\int_0^1\left(\frac{|u|^2}{2\epsilon_2}-C\right)dt.\nn\\
&\ge&\left(\frac{1}{2\epsilon_2}-\frac{1}{4\pi}\right)\|u\|^2-C\nn\\
&\ge& C_4\|u\|^2-C \lb{2.15}
\eea
for some constant $C_4>0$, where in the first inequality, we have
used (iii) of Proposition 2.5. Hence the proposition holds. \hfill\hb

By (\ref{2.15}) and the proof of Lemma 5.2.8 of \cite{Eke3}, we have

{\bf Proposition 2.8.} {\it The functional $\Psi_a$ satisfies the
Palais-Smale condition on $L_0^2(S^1, \R^{2n})$.}\hfill\hb

{\bf Proposition 2.9.} {\it $\Psi_a(u_a)<0$ for every critical
point $u_a\not= 0$ of $\Psi_a$. }

{\bf Proof}. By Proposition 2.4, we have $u_a=\dot{x}_a$ and
$x_a=\rho_ay(\tau t)$ with
\be  \frac{\varphi_a'(\rho_a)}{\rho_a}=\frac{\tau}{a}.  \lb{2.16}\ee
Hence we have
\bea \Psi_a(u_a)
&=& \int_0^1\left(\frac{1}{2}J\dot x_a\cdot x_a+G_a(-J\dot x_a)\right)dt\nn\\
&=& -\frac{1}{2}\langle H_a'(x_a),\; x_a\rangle+\int_0^1G_a(H_a'(x_a))dt\nn\\
&=& \frac{1}{2}a\varphi_a'(\rho_a)\rho_a-a\varphi_a(\rho_a). \lb{2.17}
\eea
Here the second equality follows from (\ref{2.3}) and the third
equality follows from (i) of Proposition 2.5 and (\ref{2.8}).

Let $f(t)=\frac{1}{2}a\varphi_a^\prime(t)t-a\varphi_a(t)$ for $t\ge 0$. Then
we have $f(0)=0$ and
$f'(t)=\frac{a}{2}(\varphi_a^{\prime\prime}(t)t-\varphi_a^\prime(t))<0$
since $\frac{d}{dt}(\frac{\varphi_a^\prime(t)}{t})<0$ by (iii) of
Proposition 2.2. This together with (\ref{2.16}) yield the proposition. \hfill\hb

\setcounter{equation}{0}%\setcounter{figure}{0}
\section{Critical modules for closed characteristics}%{Section 3}

In this section, we define the critical modules of  closed
characteristics and study some properties of them.

\subsection{Basic properties of critical modules}%Subsection 3.1

We have a natural $S^1$-orthogonal action on $L_0^2(S^1, \R^{2n})$ defined by
\bea\theta\cdot u(t)=u(\theta+t),\quad\forall\theta\in S^1, t\in\R.\lb{3.1}\eea
Clearly $\Psi_a$ is $S^1$-invariant.
For any $\kappa\in\R$, we denote by
\bea \Lambda_a^\kappa=\{u\in L_0^2(S^1, \R^{2n})\;|\;\Psi_a(u)\le\kappa\}.\lb{3.2}\eea
For a critical point $u$ of $\Psi_a$, we denote by
\bea \Lambda_a(u)=\Lambda_a^{\Psi_a(u)}
=\{w\in L_0^2(S^1, \R^{2n}) \;|\; \Psi_a(w)\le\Psi_a(u)\}.\lb{3.3}\eea
Clearly, both sets are $S^1$-invariant. Since the
$S^1$-action preserves $\Psi_a$, if $u$ is a critical
point of $\Psi_a$, then the whole orbit $S^1\cdot u$ is formed by
critical points of $\Psi_a$. Denote by $crit(\Psi_a)$ the set of
critical points of $\Psi_a$. Note that by the condition $(F)$,
(iii) of Proposition 2.4 and Proposition 2.6,
the number of critical orbits of $\Psi_a$ is finite.
Hence as usual we can make the following definition.

{\bf Definition 3.1.} {\it Suppose $u$ is a nonzero critical
point of $\Psi_a$, and $\Nn$ is an $S^1$-invariant
open neighborhood of $S^1\cdot u$ such that
$crit(\Psi_a)\cap(\Lambda_a(u)\cap \Nn)=S^1\cdot u$. Then
the $S^1$-critical modules of $S^1\cdot u$ is defined by
\bea C_{S^1,\; q}(\Psi_a, \;S^1\cdot u)
&=&H_{S^1,\; q}(\Lambda_a(u)\cap\Nn,\;
(\Lambda_a(u)\setminus S^1\cdot u)\cap\Nn)\nn\\
&\equiv&H_{q}((\Lambda_a(u)\cap\Nn)_{S^1},\;
((\Lambda_a(u)\setminus S^1\cdot u)\cap\Nn)_{S^1}),\lb{3.4}
\eea
where $H_{S^1,\;\ast}$ is the $S^1$-equivariant homology with
rational coefficients in the sense of A. Borel
(cf. Chapter IV of \cite{Bor1}).}

Note that this definition is independent of the choice of $\Nn$
by the excision property of the singular homology theory (cf.
Definition 1.7.5 of \cite{Cha1}). Recall that $X_{S^1}$ is defined
at the end of Section 1.

We have the following proposition for critical modules.

{\bf Proposition 3.2.} {\it The critical module $C_{S^1,\;
q}(\Psi_a, \;S^1\cdot u)$ is independent of the choice of $H_a$
defined in Proposition 2.4 in the sense that if $x_i$ are solutions
of (\ref{2.3}) with Hamiltonian functions $H_{a_i}(x)\equiv
a_i\varphi_{a_i}(j(x))$ for $i=1$ and $2$ respectively such that
both $x_1$ and $x_2$ correspond to the same closed characteristic
$(\tau, y)$ on $\Sigma$. Then we have
\be C_{S^1,\; q}(\Psi_{a_1}, \;S^1\cdot\dot {x}_1) \cong
  C_{S^1,\; q}(\Psi_{a_2}, \;S^1\cdot \dot {x}_2), \quad \forall q\in \Z.
\lb{3.5}\ee
In other words, the critical modules are invariant for all $a>\tau$ and
$\varphi_a$ satisfying (i)-(iv) of Proposition 2.2. }

{\bf Proof.} Fix a closed characteristic $(\tau, y)$ on $\Sigma$.
We assume first $\tau<a_1<a_2$. Let $\varphi_a$ be a family of functions
satisfying (i)-(iv) of Proposition 2.2 and $H_a(x)=a\varphi_a(j(x))$
satisfying Proposition 2.4 parametrized by $a\in[a_1, a_2]$. Without loss
of generality we can assume $\varphi_a$ depends continuously on $a$ in the
sense of Remark 2.3. For each $a\in[a_1, a_2]$, we denote by $x_a$ the
corresponding solution of (\ref{2.3}) with Hamiltonian $H_a$. Firstly
we prove the following

{\bf Claim. } For each $a\in[a_1, a_2]$ and $\varepsilon$ near $0$,
we have \bea |G_{a+\varepsilon}(y)-G_a(y)|&=&
O(\varepsilon)+O(\varepsilon)|y|^2,\quad \forall y\in\R^{2n},\lb{3.6}\\
|G^\prime_{a+\varepsilon}(y)-G^\prime_a(y)|&=&
O(\varepsilon)+O(\varepsilon)|y|, \quad \forall y\in\R^{2n},\lb{3.7}
\eea where we denote by $B=O(\varepsilon)$ if $|B|<C|\varepsilon|$
for some constant $C>0$.

In fact, fix an $a\in[a_1, a_2]$ and let $b\in (a-\varepsilon,
a+\varepsilon)$. For any $y\in\R^{2n}$, we have $y=\lambda
j^\prime(\xi)$ for some $\lambda\ge 0$ and $\xi\in\Sigma$. Let
$x=G_b^\prime(y)$, then by (i) of Proposition 2.5, we have $\lambda
j^\prime(\xi)=y=H_b^\prime(x)=b\varphi_b^\prime(j(x))j^\prime(x)$.
Hence $x=\mu\xi$ for some $\mu>0$. This yields
\be \lambda
j^\prime(\xi)=b\varphi_b^\prime(j(x))j^\prime(\xi).\lb{3.8} \ee
Hence $\lambda=b\varphi_b^\prime(j(x))$. Then
$j(x)=(\varphi_b^\prime)^{-1}(\lambda/b)$. Because $x=j(x)\xi$, we
obtain
\be x=(\varphi_b^\prime)^{-1}(\lambda/b)\xi, \quad
G_b^\prime(y)=x=(\varphi_b^\prime)^{-1}(\lambda/b)\xi.\lb{3.9}\ee
Hence we have
$$ |G^\prime_{a+\varepsilon}(y)-G^\prime_a(y)|=
|(\varphi_{a+\varepsilon}^\prime)^{-1}(\lambda/(a+\varepsilon))
-(\varphi_{a}^\prime)^{-1}(\lambda/a)||\xi|. $$
It suffices to consider large $|y|$, where
$(\varphi_{b}^\prime)^{-1}$ is a polynomial of degree
$1$ and  whose coefficients depend continuously on $b$
by (ii) of Proposition 2.2. Hence (\ref{3.7}) holds.

For (\ref{3.6}), we have
$$ G_b(y)=x\cdot y-H_b(x)=\lambda(\varphi_b^\prime)^{-1}(\lambda/b)
   -b\varphi_b((\varphi_b^\prime)^{-1}(\lambda/b)). $$
As above for large $|y|$ by (ii) of Proposition 2.2 we may assume
$\varphi_b$ is a polynomial of degree $2$ and whose coefficients
depend continuously on $b$, this proves (\ref{3.6}) and
the whole claim.

Now we have the following estimates:
\bea |\Psi_{a+\varepsilon}(u)-\Psi_{a}(u)|
&\le& \int_0^1|G_{a+\varepsilon}(-Ju)-G_{a}(-Ju)|dt
=O(\varepsilon)+O(\varepsilon)\|u\|^2,\lb{3.10}\\
\|\Psi_{a+\varepsilon}^\prime(u)-\Psi_{a}^\prime(u)\|^2
&=&\|JG^\prime_{a+\varepsilon}(u)-JG^\prime_a(u)\|^2
=O(\varepsilon)+O(\varepsilon)\|u\|^2.\lb{3.11} \eea
In particular, (\ref{3.10}) and (\ref{3.11}) imply that $b\mapsto\Psi_{b}$
is continuous in the $C^1$ topology. Note that the number of critical
orbits of each $\Psi_b$ is finite. Hence by the continuity of
critical modules (cf. Theorem 8.8 of \cite{MaW1} or Theorem 1.5.6 on
p.53 of \cite{Cha1}, which can be easily generalized to the equivariant
sense), our proposition holds. Note that a similar argument as above
shows that the critical modules are independent of the choice of
$\varphi_a$ in $H_a(x)=a\varphi_a(j(x))$ whenever $a$ is fixed and
$\varphi_a$ satisfies (i) to (iv) of Proposition 2.2. \hfill\hb

We say that $\Psi_a$ with $a\in [a_1, a_2]$ form a continuous family of
functionals in the sense of Proposition 3.2.

In order to compute the critical modules, as in p.35 of \cite{Eke1}
and p.219 of \cite{Eke3} we introduce the following.

{\bf Definition 3.3.} {\it Suppose $u$ is a nonzero critical
point of $\Psi_a$. Then the formal Hessian of $\Psi_a$ at $u$ is defined by
\be Q_a(v,\; v)=\int_0^1 (Jv\cdot Mv+G_a^{\prime\prime}(-Ju)Jv\cdot Jv)dt,
    \lb{3.12}\ee
which defines an orthogonal splitting $L_0^2=E_-\oplus E_0\oplus E_+$ of
$L_0^2(S^1, \R^{2n})$ into negative, zero and positive subspaces. The index
of $u$ is defined by $i(u)=\dim E_-$ and the nullity of $u$ is defined by
$\nu(u)=\dim E_0$. }

Next we show that the index and nullity defined as above are the Morse
index and nullity of a corresponding functional on a finite dimensional
subspace of $L_0^2(S^1, \R^{2n})$.

{\bf Lemma 3.4.} {\it Let  $\Psi_a$ with $a\in [a_1, a_2]$
be a continuous family of functionals defined by (\ref{2.11}).
Then there exist a finite dimensional $S^1$-invariant subspace $X$  of
$L_0^2(S^1, \R^{2n})$ and a family of $S^1$-equivariant maps
$h_a: X\rightarrow X^\perp$ such that the following hold

(i) For $g\in X$, each function $h\mapsto\Psi_a(g+h)$
has $h_a(g)$ as the unique minimum in $X^\perp$.

Let $\psi_a(g)=\Psi_a(g+h_a(g))$. Then we have

(ii) Each $\psi_a$ is $C^1$ on $X$ and $S^1$-invariant.
$g_a$ is a critical point of $\psi_a$ if and only if $g_a+h_a(g_a)$
is a critical point of $\Psi_a$.

(iii) If $g_a\in X$ and $H_a$ is $C^k$ with $k\ge 2$ in a neighborhood
of the trajectory of $g_a+h_a(g_a)$, then  $\psi_a$ is $C^{k-1}$ in a
neighborhood of $g_a$. In particular,
if $g_a$ is a nonzero critical point of $\psi_a$,
then  $\psi_a$ is $C^2$ in a neighborhood of the critical
orbit $S^1\cdot g_a$.
The index and nullity of $\Psi_a$ at $g_a+h_a(g_a)$ defined
in Definition 3.3 coincide with the Morse index and nullity of
$\psi_a$ at $g_a$.

(iv) For any $\kappa\in\R$, we denote by
\bea \widetilde{\Lambda}_a^\kappa=\{g\in X \;|\; \psi_a(g)\le\kappa\}. \lb{3.13}\eea
Then the natural embedding
$\widetilde{\Lambda}_a^\kappa \hookrightarrow  {\Lambda}_a^\kappa $
given by $g\mapsto g+h_a(g)$
is an $S^1$-equivariant homotopy equivalence.

(v) The functionals $a\mapsto\psi_a$ is continuous in $a$ in the $C^1$
topology. Moreover $a\mapsto\psi^{\prime\prime}_a$ is continuous
in a  neighborhood of the critical orbit $S^1\cdot g_a$.}

{\bf Proof.} By (ii) of Proposition 2.5, we have
\bea (G_a^\prime(u)-G_a^\prime(v), u-v)\ge
\omega|u-v|^2,\quad \forall a\in[a_1, a_2],\; u, v\in\R^{2n},
\lb{3.14}\eea
for some $\omega>0$. Hence we can use the proof of
Proposition 3.9 of \cite{Vit1} to obtain $X$ and $h_a$.
In fact, $X$ is the subspace of $L_0^2(S^1, \R^{2n})$ generated
by the eigenvectors of $-JM$ whose eigenvalues are less than
$-\frac{\omega}{2}$ and $h_a(g)$ is defined by the equation
\bea \frac{\partial}{\partial h}\Psi_a(g+h_a(g))=0,\lb{3.15}\eea
then (i)-(iii) follows from Proposition 3.9 of \cite{Vit1}.
(iv) follows from Lemma 5.1 of \cite{Vit1}.

We prove (v). As in \cite{Vit1}, (\ref{3.14}) and the
definition of $X$ yields
\bea \langle\Psi^\prime_a(u)-\Psi^\prime_a(v),\; u-v\rangle\ge
\frac{\omega}{2}\|u-v\|^2,
\quad\forall u-v\in X^\perp,\; a\in[a_1, a_2].\lb{3.16}
\eea
Hence we have
\bea  \frac{\omega}{2}\|h_{a+\epsilon}(g)-h_a(g)\|^2
&\le&\langle\Psi^\prime_{a+\epsilon}(g+h_{a+\epsilon}(g))-
\Psi^\prime_{a+\epsilon}(g+h_a(g)) ,\; h_{a+\epsilon}(g)-h_a(g)\rangle\nn\\
&=&\langle\Psi^\prime_a(g+h_a(g))-
\Psi^\prime_{a+\epsilon}(g+h_a(g)) ,\; h_{a+\epsilon}(g)-h_a(g)\rangle\nn\\
&\le&(O(\epsilon)+O(\epsilon)\|g+h_a(g)\|^2)^\frac{1}{2}\| h_{a+\epsilon}(g)-h_a(g)\|.\nn
\eea
The second equality follows by (\ref{3.15}) and the last
inequality follows by (\ref{3.11}). Hence $a\mapsto h_a(g)$
is continuous. We have $\psi_a(g)=\Psi_a(g+h_a(g))$ by definition,
$\psi_a^\prime(g)=\frac{\partial}{\partial g}\Psi_a(g+h_a(g))$
by (\ref{3.15}). Hence the first statement of (v) follows from
(\ref{3.10}) and (\ref{3.11}). The last statement of (v) follows
from p.629 of \cite{Vit1} and the implicit functional theorem
with parameters.\hfill\hb

Note that $\Psi_a$ is not $C^2$ in general, and then we can not apply
Morse theory to $\Psi_a$ directly. After the finite dimensional
approximation, the function $\psi_a$ has much better differentiability,
which allows us to apply the Morse theory to study its property.
Note that the above Lemma 3.4 is used only in Proposition 3.5 and
Theorem 4.2 below.

{\bf Proposition 3.5.} {\it For all $b\ge a>\tau$, let $\Psi_b$ be a
functional defined by (\ref{2.11}), and $u_b=\dot{x}_b$ be the critical
point of $\Psi_b$ so that $x_b$ corresponds to a fixed closed
characteristic $(\tau,y)$ on $\Sigma$ for all $b\ge a$. Then the index
$i(u_b)$ and the nullity $\nu(u_b)$ of the functional $\Psi_b$ at its
critical point $u_b$ are constants for all $b\ge a$. In particular, when
$H_b$ is $\aa$-homogenous for some $\aa\in (1,2)$ near the image set of $x_b$,
the index and nullity coincide with those defined by I. Ekeland in
\cite{Eke1} to \cite{Eke3}. Specially $1\le \nu(u_b)\le 2n-1$ always holds. }

{\bf Proof.} We consider the nullity first. As in Proposition 3.6 of
\cite{Eke1}, we have that $v\in L_0^2(S^1, \R^{2n})$ belongs to the null
space of $Q_a$ if and only if $z=Mv-J\xi$ is a solution of the linearized
system
\be \dot z(t)=JH_a^{\prime\prime}(x_a(t))z(t), \lb{3.17}\ee
for some unique $\xi\in\R^{2n}$. Denote by $R(t)$ the fundamental solution
of (\ref{3.17}). Then by Lemma 1.6.11 of \cite{Eke3}, we have
\be R(t)T_{y(0)}\Sigma\subset T_{y(\tau t)}\Sigma.\lb{3.18}\ee
Clearly, we have
\be R(1)\dot x_a(0)=\dot x_a(0).\lb{3.19}\ee
Let
\be x_a(\rho, t)=\rho y\left(\frac{\tau t}{T_\rho}\right) \quad {\rm with}\;
\frac{\tau}{T_\rho}=\frac{a\varphi_a^\prime(\rho)}{\rho}.\lb{3.20}\ee
Then we have $x_a(\rho, T_\rho)=x_a(\rho, 0)$. Differentiating
it with respect to $\rho$ and using (\ref{3.20}), we get
$$
\frac{\tau}{a}\frac{d}{d\rho}\left(\frac{\rho}{\varphi_a^\prime(\rho)}\right)\dot x_a(0)
+R(1)\rho^{-1}x_a(0)=\rho^{-1}x_a(0).$$
Hence we have
\be R(1)x_a(0)=x_a(0)-
\frac{\rho\tau}{a}\frac{d}{d\rho}\left(\frac{\rho}{\varphi_a^\prime(\rho)}\right)\dot x_a(0)
=x_a(0)+\gamma \dot x_a(0),\lb{3.21}\ee
where $\gamma<0$ since $\frac{d}{d\rho}\left(\frac{\rho}{\varphi_a^\prime(\rho)}\right)>0$
by (iii) of Proposition 2.2. For any $w\in\R^{2n}$, we have
\bea H_a^{\prime\prime}(x_a)w
&=& a\varphi_a^{\prime\prime}(j(x_a))(j^\prime(x), w)j^\prime(x)
+a\varphi_a^{\prime}(j(x_a))j^{\prime\prime}(x_a)w\nn\\
&=& a\varphi_a^{\prime\prime}(j(x_a))(j^\prime(y), w)j^\prime(y)
+\tau j^{\prime\prime}(y)w. \lb{3.22}
\eea
The last equality follows from (iii) of Proposition 2.4.
Let $z(t)=R(t)z(0)$ for $z(0)\in T_{y(0)}\Sigma$.
Then by (\ref{3.18}), we have $\dot z(t)=\tau j^{\prime\prime}(y(t))z(t)$.
Therefore $R(1)|_{T_{y(0)}}\Sigma$ is
independent of the choice of $H_a$ in Proposition 2.4. Summing up,
we have proved that in an appropriate coordinates there holds
\bea R(1)=\left(\matrix{A\quad B  \cr
                        0\quad C   \cr }\right)
\quad {\rm with}\quad A=\left(\matrix{1\quad \gamma  \cr
                        0\quad 1   \cr }\right),\nn
\eea
and $C$ is independent of $H_a$. This proves that $\nu(\dot x_a)$ is constant
for all $H_a$ satisfying Proposition 2.4 with $a>\tau$.

For any $b>a>\tau$, by (v) of Proposition 2.2, we can construct a continuous
family of $\Psi_c$ with $c\in[a, b]$ such that $H_b$ is homogenous of degree
$\aa=\aa_b$ near the image set of $x_b$. Now we can use Lemma 3.4 to obtain
a continuous family of $\psi_c$ such that $\psi_c^{\prime\prime}(g_c)$ depends
continuously on $c\in [a,b]$, where $g_c$ is the critical point of $\psi_c$
corresponding to $\dot{x}_c$. Because
$\dim\ker\psi_c^{\prime\prime}(g_c)=\nu(\dot{x_c})=constant$, the index of
$\psi_c^{\prime\prime}(g_c)$ must be constant too. Because $i(\dot{x}_b)$
and $\nu(\dot{x}_b)$ coincide with the index and nullity defined by
I. Ekeland (cf. (24) in p.219 of \cite{Eke3}), our proposition holds. \hfill\hb

In the following of this section, we fix a $\Psi_a$. All the
constructions below depend on this $\Psi_a$. In order to simplify
notations, we shall omit the subscript $a$.

In order to relate the critical modules with the index
and nullity of the critical point, we use the finite
dimensional approximation introduced by I. Ekeland in \cite{Eke1}.
For $\epsilon>0$, we define $\Psi_{a, \epsilon}(v)\equiv
\int_0^\epsilon(\frac{1}{2}Jv\cdot M_\epsilon v+G_a(-Jv))dt$
for $v\in L^2([0, \epsilon], \R^{2n})$, where
$M_\epsilon v(t)=\int_0^t v(s)ds$. Then we have

{\bf Proposition 3.6.}(cf. Lemma 3.9 of \cite{Eke1})
{\it For $\xi\in\R^{2n}$ and $\epsilon>0$, consider the problem
\be
\min\left\{\Psi_{a, \epsilon}(v)\;\left|\frac{}{}\right.\;\int_0^\epsilon v(t)dt=\xi\right\}.
\lb{3.23}\ee
Then there exists $\epsilon_0>0$ such that for all $\epsilon<\epsilon_0$
and $\xi\in\R^{2n}$, the problem (\ref{3.23}) has a unique
solution $v(\epsilon, \xi)$ which is $C^2$ on $t$.
We have $v(\epsilon, 0)=0$ for all $\epsilon$ and if $\xi\neq0$,
then $v(\epsilon, \xi)(t)\neq 0$ for all $t$. }\hfill\hb

Now following I. Ekeland, we choose an $\iota\in\N$ such that $\iota^{-1}<\epsilon_0$
and let $\epsilon=\iota^{-1}$. Define
\be (\R^{2n})_0^\iota=
\left\{(\xi_1, \ldots, \xi_\iota)\in(\R^{2n})^\iota\;\left|\frac{}{}\right.\;
\sum_{i=1}^\iota \xi_i=0\right\}.\lb{3.24}\ee
Let $\Omega=\{(\xi_1, \ldots, \xi_\iota)\in(\R^{2n})_0^\iota\;|\; \xi_i\neq 0 ,\;\forall i\}$.
Let $p_\iota: L_0^2(S^1, \R^{2n})\rightarrow(\R^{2n})_0^\iota$ to be
\be p_\iota v=\left(\int_0^\epsilon vdt,
\int_\epsilon^{2\epsilon}vdt,\ldots,\int_{(\iota-1)\epsilon}^1 vdt\right).\lb{3.25}\ee
Let $r_\iota: (\R^{2n})_0^\iota\rightarrow L_0^2(S^1, \R^{2n})$ to be
\be r_\iota(\xi_1, \ldots, \xi_\iota)(t)=
v(\epsilon, \xi_k)(t-(k-1)\epsilon),\quad (k-1)\epsilon\le t\le k\epsilon,\lb{3.26}
\ee
where $v(\epsilon, \xi_k)$ is given by Proposition 3.6.

{\bf Lemma 3.7.} (cf. Lemma 3.10 and 3.11 of \cite{Eke1}) {\it
The functional $\Psi_a\circ r_\iota$ is $C^1$ on $(\R^{2n})_0^\iota$
and $C^2$ on $\Omega$. If $(\xi_1, \ldots, \xi_\iota)\neq 0$
is a critical point of $\Psi_a\circ r_\iota$, then
$r_\iota(\xi_1, \ldots, \xi_\iota)\neq 0$ is a critical point of
$\Psi_a$. Conversely, if $u\neq 0$ is a critical point of
$\Psi_a$, then $p_\iota u$ belongs to $\Omega$ and is a critical
point of $\Psi_a\circ r_\iota$. Moreover, $u$ and $p_\iota u$ have
the same index and nullity.  }\hfill\hb

Now let $u\neq 0$ be a critical point of $\Psi_a$
with multiplicity $mul(u)=m$, i.e., $u$ corresponds
to a closed characteristic $(m\tau, y)\subset\Sigma$
with $(\tau, y)$ being  prime. Hence $u(t+\frac{1}{m})=u(t)$
for all $t\in \R$ and the orbit of $u$, namely,
$S^1\cdot u\cong S^1/\Z_m\cong S^1$.
Let $p: N(S^1\cdot u)\rightarrow S^1\cdot u$ be the normal
bundle of $S^1\cdot u$ in $L_0^2(S^1, \R^{2n})$ and let
$p^{-1}(\theta\cdot u)=N(\theta\cdot u)$ be the fibre over
$\theta\cdot u$, where $\theta\in S^1$.
Let $DN(S^1\cdot u)$ be the $\varrho$ disk bundle of $N(S^1\cdot u)$
for some $\varrho>0$ sufficiently small, i.e.,
$DN(S^1\cdot u)=\{\xi\in N(S^1\cdot u)\;| \; \|\xi\|<\varrho\}$
which is identified by the exponential map with a subset of $L_0^2(S^1, \R^{2n})$,
and let $DN(\theta\cdot u)=p^{-1}(\theta\cdot u)\cap DN(S^1\cdot u)$
be the disk over $\theta\cdot u$. Clearly, $DN(\theta\cdot u)$ is
$\Z_m$-invariant and we have $DN(S^1\cdot u)=DN(u)\times_{\Z_m}S^1$
where the $Z_m$ action is given by
$$(\theta, v, t)\in \Z_m\times DN(u)\times S^1\mapsto
(\theta\cdot v, \;\theta^{-1}t)\in DN(u)\times S^1.$$
Hence for an $S^1$ invariant subset $\Gamma$ of $DN(S^1\cdot u)$,
we have $\Gamma/S^1=(\Gamma_u\times_{\Z_m}S^1)/S^1=\Gamma_u/\Z_m$,
where $\Gamma_u=\Gamma\cap DN(u)$. Let $\Gamma(\iota)=\im r_\iota$
and for any $\kappa\in\R$, we denote by
\bea \Gamma(\iota)^\kappa=\{u\in \Gamma(\iota)\;|\; \Psi_a(u)\le\kappa\}. \lb{3.27}\eea

{\bf Lemma 3.8.} {\it For any $\kappa\in\R$, we have a
$\Z_\iota$-equivariant deformation retract from
$\Lambda_a^\kappa$ to $\Gamma(\iota)^\kappa$.}

{\bf Proof.}  For any $v\in L_0^2(S^1, \R^{2n})$,
let $p_\iota(v)=(\xi_1, \ldots, \xi_\iota)$, we have
\bea \Psi_a(v)&=&\int_0^1\left(\frac{1}{2}Jv\cdot Mv+G_a(-Jv)\right)dt\nn\\
&=&\sum_{j=1}^\iota\int_{(j-1)\epsilon}^{j\epsilon}\left(
\frac{1}{2}Jv\cdot\left(\int_0^{(j-1)\epsilon}v(s)ds+
\int_{(j-1)\epsilon}^tv(s)ds\right )+G_a(-Jv)\right)dt\nn\\
&=&\sum_{j=1}^\iota\left(\Psi_{a, \epsilon}(v_j)
+\frac{1}{2}J\xi_j\cdot\sum_{l=1}^{j-1}\xi_l\right), \lb{3.28}
\eea
where $v_j(t)=v(t+(j-1)\epsilon)$.
By Proposition 3.6,  we have $\Psi_a(v)\ge \Psi_a(r_\iota p_\iota v)$.
Hence the deformation retract
$F: \Lambda_a^\kappa\times[0, 1]\rightarrow\Lambda_a^\kappa $
is given by $(v, s)\mapsto sr_\iota p_\iota v+(1-s)v$. This is well defined,
since by Lemma 3.9 of \cite{Eke1}, we have $\Psi_{a, \epsilon}$
is strictly convex, hence
$$\Psi_a(sr_\iota p_\iota v+(1-s)v)
\le s\Psi_a(r_\iota p_\iota v)+(1-s)\Psi_a(v)\le \Psi_a(v).$$
Clearly $F$ is $\Z_\iota$-equivariant and $F=id$ on $\Gamma(\iota)^\kappa$,
hence the lemma holds.\hfill\hb

As in p.51 of \cite{Rad2}, let
\bea D_\iota N(S^1\cdot u)=DN(S^1\cdot u)\cap\Gamma(\iota), \quad
D_\iota N(\theta\cdot u)=DN(\theta\cdot u)\cap\Gamma(\iota).\lb{3.29}
\eea
For a $\Z_m$-space pair $(A, B)$, let
\be H_{\ast}(A, B)^{\pm\Z_m}
= \{\sigma\in H_{\ast}(A, B)\,|\,L_{\ast}\sigma=\pm \sigma\}, \lb{3.30}\ee
where $L$ is a generator of the $\Z_m$-action. Then by Lemma 3.8,
as in Section 6 of \cite{Rad2} or Section 3 of \cite{BaL1}, we have

{\bf Lemma 3.9.} {\it Suppose $u\neq 0$ is a critical point
of $\Psi_a$ with $mul(u)=m$ and suppose $\Gamma(\iota)$ is a
finite dimensional approximation as above with $m|\iota$, i.e.,
$m$ is a factor of $\iota$. Then we have}
\bea C_{S^1,\; \ast}(\Psi_a, \;S^1\cdot u)
&\cong& H_\ast((\Lambda_a(u)\cap DN(u))/\Z_m,\;
    ((\Lambda_a(u)\setminus\{u\})\cap DN(u))/\Z_m) \nn\\
&\cong& H_\ast((\Lambda_a(u)\cap D_\iota N(u))/\Z_m,\;
    ((\Lambda_a(u)\setminus\{u\})\cap D_\iota N(u))/\Z_m) \nn\\
&\cong& H_\ast(\Lambda_a(u)\cap D_\iota N(u),\;
    (\Lambda_a(u)\setminus\{u\})\cap D_\iota N(u))^{\Z_m}.\lb{3.31}\eea

{\bf Proof.} For reader's conveniences, we sketch a proof here and refer to
Section 6 of \cite{Rad2} or Section 3 of \cite{BaL1} for related details.

By Definition 3.1, we have
\bea C_{S^1,\; \ast}(\Psi_a, \;S^1\cdot u)\cong
H_{S^1, \;\ast}(\Lambda_a(u)\cap DN(S^1\cdot u),\;
(\Lambda_a(u)\setminus\{S^1\cdot u\})\cap DN(S^1\cdot u)). \nn
\eea
Since all the isotropy groups $G_x=\{g\in S^1\;|\;g\cdot x=x\}$
for $x\in DN(S^1\cdot u)$ are finite, we can use Lemma 6.11 of
\cite{FaR1} to obtain
\bea
&& H_{S^1}^\ast(\Lambda_a(u)\cap DN(S^1\cdot u),\;
(\Lambda_a(u)\setminus\{S^1\cdot u\})\cap DN(S^1\cdot u))\nn\\
&&\qquad \cong H^\ast ((\Lambda_a(u)\cap DN(S^1\cdot u))/S^1, \;
((\Lambda_a(u)\setminus\{S^1\cdot u\})\cap DN(S^1\cdot u))/S^1)\nn\\
&&\qquad \cong H^\ast ((\Lambda_a(u)\cap DN(u))/\Z_m, \;
((\Lambda_a(u)\setminus\{u\})\cap DN(u))/\Z_m).  \nn \eea
By the condition (F) at the beginning of Section 2, a small perturbation
on the energy functional can be applied to reduce each critical
orbit to nearby non-degenerate ones. Thus similar to the proofs of Lemma
2 of \cite{GrM1} and Lemma 4 of \cite{GrM2}, all the homological
$\Q$-modules of each space pair in the above relations are all finitely
generated. Therefore we can apply Theorem 5.5.3 and Corollary 5.5.4 on
pages 243-244 of \cite{Spa1} to obtain the same relation on homological
$\Q$-modules:
\bea
&& H_{S^1,\ast}(\Lambda_a(u)\cap DN(S^1\cdot u),\;
(\Lambda_a(u)\setminus\{S^1\cdot u\})\cap DN(S^1\cdot u))\nn\\
&&\qquad \cong H_{\ast}((\Lambda_a(u)\cap DN(S^1\cdot u))/S^1, \;
((\Lambda_a(u)\setminus\{S^1\cdot u\})\cap DN(S^1\cdot u))/S^1)\nn\\
&&\qquad \cong H_{\ast}((\Lambda_a(u)\cap DN(u))/\Z_m, \;
((\Lambda_a(u)\setminus\{u\})\cap DN(u))/\Z_m).   \nn
\eea
Now by Lemma 3.8, as in Section 6.2 of \cite{Rad2} or Section 3
of \cite{BaL1}, we obtain
\bea
&& H_\ast((\Lambda_a(u)\cap DN(u))/\Z_m,\;
    ((\Lambda_a(u)\setminus\{u\})\cap DN(u))/\Z_m) \nn\\
&&\qquad \cong H_\ast((\Lambda_a(u)\cap D_\iota N(u))/\Z_m,\;
    ((\Lambda_a(u)\setminus\{u\})\cap D_\iota N(u))/\Z_m). \nn\eea
Note that the same argument as in Section 6.3 of \cite{Rad2},
in particular Satz 6.6 of \cite{Rad2}, Lemma 3.6 of \cite{BaL1} or
Theorem 3.2.4 of \cite{Bre1} yields
\bea
&& H_\ast((\Lambda_a(u)\cap D_\iota N(u))/\Z_m,\;
    ((\Lambda_a(u)\setminus\{u\})\cap D_\iota N(u))/\Z_m) \nn\\
&&\qquad \cong H_\ast(\Lambda_a(u)\cap D_\iota N(u),\;
    (\Lambda_a(u)\setminus\{u\})\cap D_\iota N(u))^{\Z_m}. \nn\eea
The above relations together complete the proof of Lemma 3.9.
\hfill\hb

\subsection{The periodic property of critical modules}%Subsection 3.2

By Proposition 3.2, we have that $C_{S^1,\; \ast}(\Psi_a, \;S^1\cdot u)$
is independent of the choice of the Hamiltonian function $H_a$ whenever
$H_a$ satisfies conditions in Proposition 2.4. Hence in order to compute
the critical modules, we can choose $\Psi_a$ with $H_a$ being positively
homogeneous of degree $\aa=\aa_a$ near the image set of every nonzero solution
$x$ of (\ref{2.3}) which corresponding to some closed characteristic $(\tau,y)$
with period $\tau$ being strictly less than $a$.

In other words, for a given $a>0$, we choose $\vartheta\in (0,1)$ first such that
$[a\vartheta, a(1-\vartheta)]\supset per(\Sigma)\cap (0,a)$ holds by the
definition of the set $per(\Sg)$ and the assumption (F). Then we choose
$\aa=\aa_a\in (1,2)$ sufficiently close to $2$ by (v) of Proposition 2.2 such
that $\varphi_a(t)=ct^\alpha$ for some constant $c>0$ and $\alpha\in(1,\,2)$
whenever $\frac{\varphi_a^\prime(t)}{t}\in [\vartheta, 1-\vartheta]$.
In this subsection we suppose that $\varphi_a$ possesses this property (v).

Now we consider iterations of critical points of $\Psi_a$.
Suppose $u\neq 0$ is a critical point of $\Psi_a$ with $mul(u)=m$.
By Propositions 2.4 and 2.6, we have $u=\dot x$ with $x$ being a
solution of (\ref{2.3}) and $x=\rho y(\tau t)$ with
$\frac{\varphi_a^\prime(\rho)}{\rho}=\frac{\tau}{a}$.
Moreover, $(\tau, y)$ is a closed characteristic on $\Sigma$
with minimal period $\frac{\tau}{m}$.
For any $p\in\N$ satisfying $p\tau<a$, the $p$th iteration $u^p$ of $u$
is given by $\dot x^p$, where $x^p$ is the unique
solution of (\ref{2.3}) corresponding to $(p\tau, y)$.
Hence we have
\bea
&& x(t)=\left(\frac{\tau}{c\alpha a}\right)^\frac{1}{\alpha-2}y(\tau t), \quad
x^p(t)=\left(\frac{p\tau}{c\alpha a}\right)^\frac{1}{\alpha-2}y(p\tau t),\nn\\
&& u(t)=\dot x(t)=
\tau^\frac{\alpha-1}{\alpha-2}{(c\alpha a)}^\frac{1}{2-\alpha}\dot y(\tau t), \quad
u^p(t)=\dot x^p(t)=
(p\tau)^\frac{\alpha-1}{\alpha-2}{(c\alpha a)}^\frac{1}{2-\alpha}\dot y(p\tau t).\nn
\eea
These yield \be u^p(t)=p^\frac{\alpha-1}{\alpha-2}u(pt).\lb{3.32}\ee
Choose two finite dimensional approximations $\Gamma(\iota)$ and
$\Gamma(p\iota)$ as above and define the $p$th iteration $\phi^p$
on $D_\iota N(u)$ by
\bea \phi^p: v(t)\mapsto p^\frac{\alpha-1}{\alpha-2}v(pt).\lb{3.33}
\eea

{\bf Claim}.  $\phi^p$ maps $D_\iota N(u)$ into $D_{p\iota }N(u^p)$
if the radii of the two normal disk bundles are suitably chosen.

In fact, clearly $\phi^p(v)\in DN(u^p)$
and by Lemma 3.9 of \cite{Eke1},
\bea v\left(t+\frac{k-1}{\iota}\right)=
v\left(\frac{1}{\iota}, \int_\frac{k-1}{\iota}^\frac{k}{\iota}v(s)ds\right)(t)\equiv v_k(t),
\quad 1\le k\le \iota,\; 0\le t\le\frac{1}{\iota}.\lb{3.34}\eea
Here $v_k$ is the unique solution of
\bea v_k(t)=JH_a^\prime(M_{\frac{1}{\iota}}v_k+\zeta_k),\quad
\int_0^\frac{1}{\iota}v_k(t)dt=\int_\frac{k-1}{\iota}^\frac{k}{\iota}v(t)dt,\lb{3.35}\eea
where $\zeta_k$ is uniquely determined by $v_k$.
Hence it suffices to show that
\bea  \phi^p(v)(t)=JH_a^\prime(M_\frac{1}{p\iota}\phi^p(v)(t)+\zeta^\prime),
\quad t\in\left[\frac{l\iota+j-1}{p\iota}, \frac{l\iota+j}{p\iota}\right],  \nn
\eea
for some $\zeta^\prime\in\R^{2n}$ and $0\le l<p,\;1\le j<\iota$.
An easy computation show that
\bea M_\frac{1}{p\iota}(\phi^p(v))(t)=
p^\frac{1}{\alpha-2}(M_\frac{1}{\iota}v)(pt).\lb{3.36}\eea
Then we have
\bea \phi^p(v)(t)&=&p^\frac{\alpha-1}{\alpha-2}v(pt)\nn\\
&=&p^\frac{\alpha-1}{\alpha-2}JH_a^\prime((M_{\frac{1}{\iota}}v_j)(pt)+\zeta_j)\nn\\
&=&JH_a^\prime(p^\frac{1}{\alpha-2}((M_{\frac{1}{\iota}}v_j)(pt)+\zeta_j))\nn\\
&=&JH_a^\prime(M_\frac{1}{p\iota}\phi^p(v)(t)+\zeta^\prime).\nn
\eea
In the above computations, we have used that $H_a$ and then $H_a^\prime$ are
positively homogeneous of degrees $\alpha$ and $\alpha-1$ respectively. This is
true since by Proposition 3.6, all $v\in D_\iota N(u^i)$ lies in an $L^\infty$
neighborhood of $u^i$ for $1\le i\le p$. This proves the claim.

We have
\bea \Psi_a(\phi^p(v))&=&\int_0^1\left(
\frac{1}{2}J\phi^p(v)(t)\cdot M\phi^p(v)(t)+G_a(-J\phi^p(v)(t))\right)dt\nn\\
&=& \int_0^1\left(\frac{1}{2}Jp^\frac{\alpha-1}{\alpha-2}v(pt)\cdot
p^\frac{1}{\alpha-2}(Mv)(pt)+G_a(-Jp^\frac{\alpha-1}{\alpha-2}v(pt))\right)dt\nn\\
&=& p^\frac{\alpha}{\alpha-2}\int_0^1\left(\frac{1}{2}Jv(pt)\cdot
(Mv)(pt)+G_a(-Jv(pt))\right)dt \nn\\
&=& p^{\frac{\alpha}{\alpha-2}}\Psi_a(v). \lb{3.37}
\eea
Here the second equality follows from (\ref{3.36}) and the
third equality follows from (v) of Proposition 2.5.
We define a new inner product $\langle\cdot, \cdot\rangle_p$ on
$L^2_0(S^1, \R^{2n})$ by
\bea \langle v, w\rangle_p=p^\frac{\alpha-2}{2(\alpha-1)}\langle v, w\rangle.\lb{3.38}
\eea
Then $\phi^p: D_\iota N(u)\rightarrow D_{p\iota}N(u^p)$ is an isometry
from the standard inner product to the above one.
Clearly $\phi^p(D_\iota N(u))$ consists of points in $D_{p\iota}N(u^p)$
which are fixed by the $\Z_p$-action. Since the
$\Z_p$-action on $D_{p\iota}N(u^p)$ are isometries and
$f\equiv\Psi_a|_{D_{p\iota}N(u^p)}$
is $\Z_p$-invariant, we have
\bea f^{\prime\prime}(v)
=\left(\matrix{(f|_{\phi^p(D_\iota N(u))})^{\prime\prime} \quad 0\cr
\qquad 0\qquad\qquad\;\;\ast}\right),\quad \forall v\in \phi^p(D_\iota N(u)).
\lb{3.39}\eea
Moreover, we have
\bea f^\prime(v)=(f|_{\phi^p(D_\iota N(u))})^{\prime},\quad
\forall v\in \phi^p(D_\iota N(u)).\lb{3.40}\eea
Now we can apply the results by D. Gromoll and W. Meyer
\cite{GrM1}  to the manifold $D_{p\iota}N(u^p)$, with
$u^p$ as its unique critical point. Then
$mul(u^p)=pm$ is the multiplicity of $u^p$ and the isotropy  group
$\Z_{pm}\subseteq S^1$ of $u^p$ acts on $D_{p\iota}N(u^p)$ by isometries.
According to  Lemma 1 of \cite{GrM1},  we have a $\Z_{pm}$-invariant
decomposition of $T_{u^p}(D_{p\iota}N(u^p))$
\bea T_{u^p}(D_{p\iota}N(u^p))
=V^+\oplus V^-\oplus V^0=\{(x_+, x_-, x_0)\}\lb{3.41}
\eea
with $\dim V^-=i(u^p)$, $\dim V^0=\nu(u^p)-1$ and a
$\Z_{pm}$-invariant neighborhood $B=B_+\times B_-\times B_0$
for $0$ in $T_{u^p}(D_{p\iota}N(u^p))$ together with
two $\Z_{pm}$-invariant diffeomorphisms
$$\Phi :B=B_+\times B_-\times B_0\rightarrow
\Phi(B_+\times B_-\times B_0)\subset D_{p\iota}N(u^p)$$
and
$$ \eta : B_0\rightarrow W(u^p)\equiv\eta(B_0)\subset D_{p\iota}N(u^p)$$
and $\Phi(0)=\eta(0)=u^p$, such that
\be \Psi_a\circ\Phi(x_+,x_-,x_0)=|x_+|^2 - |x_-|^2 + \Psi_a\circ\eta(x_0),
\lb{3.42}\ee
with $d(\Psi_a\circ \eta)(0)=d^2(\Psi_a\circ\eta)(0)=0$. As usual, we
call $W(u^p)$ a local characteristic manifold and $U(u^p)=B_-$ a local
negative disk at $u^p$. By the proof of Lemma 1 of \cite{GrM1},
$W(u^p)$ and $U(u^p)$ are $\Z_{pm}$-invariant. It follows from
(\ref{3.42}) that $u^p$ is an isolated critical point of
$\Psi_a|_{D_{p\iota}N(u^p)}$. Then as in Lemma 6.4 of \cite{Rad2},
we have
\bea &&H_\ast(\Lambda_a(u^p)\cap D_{p\iota}N(u^p),\;
(\Lambda_a(u^p)\setminus\{u^p\})\cap D_{p\iota}N(u^p))\nn\\
=&&H_\ast (U(u^p),U(u^p)\setminus\{u^p\}) \otimes
H_\ast(W(u^p)\cap \Lambda_a(u^p), (W(u^p)\setminus\{u^p\})\cap \Lambda_a(u^p)),
\lb{3.43}\eea
where \be H_q(U(u^p),U(u^p)\setminus\{u^p\} )
    = \left\{\matrix{\Q, & {\rm if\;}q=i(u^p),  \cr
                      0, & {\rm otherwise}. \cr}\right.  \lb{3.44}\ee
Now we have the following proposition.

{\bf Proposition 3.10.} {\it Let $u\neq 0$ be a critical
point of $\Psi_a$ with $mul(u)=1$.  Then for all $p\in\N$
and $q\in\Z$, we have
\be C_{S^1,\; q}(\Psi_a, \;S^1\cdot u^p)\cong
\left(\frac{}{}H_{q-i(u^p)}(W(u^p)\cap \Lambda_a(u^p),
(W(u^p)\setminus\{u^p\})\cap \Lambda_a(u^p))\right)^{\beta(u^p)\Z_p},
\lb{3.45}\ee
where $\beta(u^p)=(-1)^{i(u^p)-i(u)}$. In particular, if $u^p$ is
non-degenerate, i.e., $\nu(u^p)=1$, then}
\be C_{S^1,\; q}(\Psi_a, \;S^1\cdot u^p)
    = \left\{\matrix{\Q, & {\rm if\;}q=i(u^p)\;{\rm and\;}\beta(u^p)=1,  \cr
                      0, & {\rm otherwise}. \cr}\right.  \lb{3.46}\ee

{\bf Proof.} Suppose $\theta$ is a generator of the linearized
$\Z_p$-action on $U(u^p)$. Then  $\theta(\xi)=\xi$ if and only if
$\xi\in T_{u^p}({\phi^p(D_\iota N(u))})$.
Hence it follows from (\ref{3.37}) and (\ref{3.39})
that $\xi=(\phi^p)_\ast(\xi^\prime)$
for a unique
$\xi^\prime\in T_u(D_\iota N(u))^-$.
Hence the proof of Satz 6.11 in \cite{Rad2} or Proposition 2.8
in \cite{BaL1} yield this proposition.
\hfill\hb

{\bf Definition 3.11.} {\it Let $u\neq 0$ be a critical
point of $\Psi_a$ with $mul(u)=1$.
Then for all $p\in\N$ and $l\in\Z$, let
\bea
k_{l, \pm 1}(u^p)&=&\dim\left(\frac{}{}H_l(W(u^p)\cap \Lambda_a(u^p),
(W(u^p)\setminus\{u^p\})\cap \Lambda_a(u^p))\right)^{\pm\Z_p},\lb{3.47}\\
k_l(u^p)&=&\dim\left(\frac{}{}H_l(W(u^p)\cap \Lambda_a(u^p),
(W(u^p)\setminus\{u^p\})\cap \Lambda_a(u^p))\right)^{\beta(u^p)\Z_p}.\lb{3.48}\eea
$k_l(u^p)$s are called critical type numbers of $u^p$. }

Note that by Proposition 3.5, we have $k_{l, \pm 1}(u^p)=0$ if
$l\notin [0, 2n-2]$.

Similar to Section 7.1 of \cite{Rad2} or Theorem 2.11 of \cite{BaL1},
we have

{\bf Lemma 3.12.} {\it Let $u\neq 0$ be a critical point of $\Psi_a$
with $mul(u)=1$. Suppose $\nu(u^m)=\nu(u^{pm})$ for some $m, p\in\N$,
then we have $k_{l, \pm 1}(u^m)=k_{l, \pm 1}(u^{pm})$ for all $l\in\Z$.}

{\bf Proof.} We choose finite dimensional approximations
$\Gamma(\iota)$ and $\Gamma(p\iota)$ as in Proposition 3.6 with
$m|\iota$ and let $\phi^p: D_\iota N(u^m)\rightarrow
D_{p\iota}N(u^{pm})$ be the $p$th iteration map. By (\ref{3.38}),
$\phi^p$ is an isometry under the modified metric. Hence by
(\ref{3.37}),  we have
\be \nu(u^m)-1=\dim\ker((\Psi_a|_{D_\iota
N(u^m)})^{\prime\prime}-I) =\dim\ker((\Psi_a|_{\phi^p(D_\iota
N(u^m))})^{\prime\prime}-I). \lb{3.49}\ee
Thus by (\ref{3.39}) and the assumption $\nu(u^m)=\nu(u^{pm})$, we have
that $T_{u^{pm}}(\phi^p(D_\iota N(u^m)))$ contains the nullity space of
the Hessian of $\Psi_a|_{D_{p\iota}N(u^{pm})}$. Now by
(\ref{3.40}), we can use Lemma 7 of \cite{GrM1} to obtain that
$\phi^p(W(u^m))\equiv W(u^{pm})$ is a characteristic manifold of
$\Psi_a|_{D_{p\iota}N(u^{pm})}$, where $W(u^m)$ is a characteristic
manifold of $\Psi_a|_{D_{\iota}N(u^m)}$. By (\ref{3.37}), we have
\bea \phi^p:
&&(W(u^m)\cap \Lambda_a(u^m),(W(u^m)\setminus\{u^m\})\cap \Lambda_a(u^m)) \nn\\
&&\qquad \rightarrow (W(u^{pm})\cap \Lambda_a(u^{pm}),(W(u^{pm})
  \setminus\{u^{pm}\})\cap \Lambda_a(u^{pm})) \nn
\eea
is a homeomorphism. Suppose $\theta$ and $\theta_p$ generate
the $\Z_m$ and $\Z_{pm}$ action on $W(u^m)$ and $W(u^{pm})$
respectively. Then clearly
$\phi^p\circ \theta=\theta_p\circ \phi^p$ holds and it implies
\bea
&& H_\ast(W(u^m)\cap \Lambda_a(u^m),
 (W(u^m)\setminus\{u^m\})\cap \Lambda_a(u^m))^{\pm\Z_m}\nn\\
&&\qquad \cong (W(u^{pm})\cap \Lambda_a(u^{pm}),
 (W(u^{pm})\setminus\{u^{pm}\})\cap \Lambda_a(u^{pm}))^{\pm\Z_{pm}}.
\nn\eea
Therefore our lemma holds. \hfill\hb

{\bf Proposition 3.13.} {\it Let $u\neq 0$ be a critical point of
$\Psi_a$ with $mul(u)=1$. Then there exists a minimal
$K(u)\in 2\N$ such that
\bea
&& \nu(u^{p+K(u)})=\nu(u^p),\quad i(u^{p+K(u)})-i(u^p)\in 2\Z,
     \qquad\forall p\in \N,  \lb{3.50}\\
&& k_l(u^{p+K(u)})=k_l(u^p), \qquad\forall p\in \N,\;l\in\Z. \lb{3.51}\eea
We call $K(u)$ the minimal period of critical modules of iterations of the
functional $\Psi_a$ at $u$. }

{\bf Proof.} As in the proof of Proposition 3.5, we denote by
$R(t)$ the fundamental solution of (\ref{3.17}). Then by
Lemma 1.1 and 1.2 of \cite{LoZ1}, we have $i(u^p)=i(u,p)-n$
and $\nu(u^p)=\nu(u,p)$ for all $p\in\N$, where
$(i(u,p), \nu(u,p))$ are index and nullity defined
by C. Conley and E. Zehnder in \cite{CoZ1}, Y. Long and
E. Zehnder in \cite{LZe1} and Y. Long in \cite{Lon1},
cf. \cite{Lon4}. Hence we have $\nu(u^p)=\dim\ker (R(1)^p-I_{2n})$.
Denote by $\lambda_i=exp(\pm2\pi\frac{r_i}{s_i})$
the eigenvalues of $R(1)$ possessing rotation angles which are rational
multiple of $\pi$ with $r_i$, $s_i\in\N$ and $(r_i,s_i)=1$ for
$1\le i\le q$. Let $K(u)$ be twice of the least common multiple of
$s_1,\ldots, s_q$. Then (\ref{3.50}) holds.
Note that the later conclusion in (\ref{3.50}) follows from Theorem 9.3.4
of \cite{Lon4}.

In order to prove (\ref{3.51}), it suffices to show
\be k_l(u^{m+qK(u)})=k_l(u^m), \qquad \forall q\in \N, \; l\in\Z,
         \; 1\le m\le K(u).  \lb{3.52}\ee

In fact, assume that (\ref{3.52}) is proved. Note that (\ref{3.51}) follows
from (\ref{3.52}) with $q=1$ directly when $p\le K(u)$. When $p>K(u)$, we
write $p=m+qK(u)$ for some $q\in\N$ and $1\le m\le K(u)$. Then by (\ref{3.52})
we obtain
$$ k_l(u^{p+K(u)})=k_l(u^{m+(q+1)K(u)})=k_l(u^m)=k_l(u^{m+qK(u)})=k_l(u^p), $$
i.e., (\ref{3.51}) holds.

To prove (\ref{3.52}), we fix an integer $m\in [1,K(u)]$. Let
$$ A=\{s_i\in\{s_1,\ldots, s_q\}\;|\;\;s_i\;{\rm is\;a\;factor\;of\;}m\},  $$
and let $m_1$ be the least common multiple of elements in $A$. Hence we have
$m=m_1m_2$ for some $m_2\in \N$ and $\nu(u^m)=\nu(u^{m_1})$. Thus by Lemma
3.12, we have $k_l(u^m)=k_{l,\beta(u^m)}(u^{m_1})$. Since $m+pK(u)=m_1m_3$
for some $m_3\in\N$, we have by Lemma 3.12 that
$k_l(u^{m+pK(u)})=k_{l, \beta(u^{m+pK(u)})}(u^{m_1})$.
By (\ref{3.50}), we obtain $\beta(u^{m+pK(u)})=\beta(u^m)$, and then
(\ref{3.52}) is proved. This completes the proof. \hfill\hb

Note that the above Proposition 3.13 could be established also without forcing
the Hamiltonian to be homogeneous near its critical points. In fact,
by Proposition 3.2, it holds for any Hamiltonian defined by
Proposition 2.4.

\subsection{Indices and Euler characteristics of closed characteristics }%Subsection 3.3

In the following, Let $\Psi_a$ by any function defined by (\ref{2.11})
with $H_a$ satisfying Proposition 2.4, we do not require $H_a$
to be homogeneous anymore.

{\bf Definition 3.14.} {\it Suppose the condition (F) at the beginning of \S 2
holds. For every closed characteristic $(\tau,y)$ on $\Sigma$, let $a>\tau$ and
choose $\vf_a$ to satisfy (i)-(iv) of Proposition 2.2.  Determine $\rho$ uniquely by
$\frac{\vf_a'(\rho)}{\rho}=\frac{\tau}{a}$. Let $x=\rho y(\tau t)$ and $u=\dot{x}$.
Then we define the index $i(\tau,y)$ and nullity $\nu(\tau,y)$ of $(\tau,y)$ by
$$ i(\tau,y)=i(u), \qquad \nu(\tau,y)=\nu(u). $$
Then the mean index of $(\tau, y)$ is defined by }
\be \hat i(\tau, y)=\lim_{m\rightarrow\infty}\frac{i(m\tau, y)}{m}.\lb{3.53}\ee

Note that by Proposition 3.5, the index and nullity are well defined and is
independent of the choice of $a>\tau$ and $\vf_a$ satisfying (i)-(iv) of
Proposition 2.2. Note that by Theorem 1.7.7 of \cite{Eke3}
(cf. Corollary 8.3.2 of \cite{Lon4}), we have $\hat i(\tau, y)>2$.

For a prime closed characteristic $(\tau, y)$ on $\Sigma$,
we denote simply by $y^m\equiv(m\tau, y)$ for $m\in\N$.
By Proposition 3.2, we can define the critical
type numbers $k_l(y^m)$ of $y^m$ to be $k_l(u^m)$, where $u^m$ is
the critical point of $\Psi_a$ corresponding to $y^m$.
We also define $K(y)=K(u)$, where $K(u)\in \N$ is given by
Proposition 3.13. Suppose $\Nn$ is an $S^1$-invariant
open neighborhood of $S^1\cdot u^m$ such that
$crit(\Psi_a)\cap(\Lambda_a(u^m)\cap\Nn)=S^1\cdot u^m$.
Then we make the following definition

{\bf Definition 3.15.} {\it The Euler characteristic $\chi(y^m)$
of $y^m$ is defined by
\bea \chi(y^m)
&\equiv& \chi((\Lambda_a(u^m)\cap\Nn)_{S^1},\;
((\Lambda_a(u^m)\setminus S^1\cdot u^m)\cap\Nn)_{S^1}) \nn\\
&\equiv& \sum_{q=0}^{\infty}(-1)^q\dim C_{S^1,\; q}(\Psi_a,\;S^1\cdot u^m).
\lb{3.54}\eea
Here $\chi(A, B)$ denotes the usual Euler characteristic of the space pair $(A, B)$.

The average Euler characteristic $\hat\chi(y)$ of $y$ is defined by }
\be \hat{\chi}(y)=\lim_{N\to\infty}\frac{1}{N}\sum_{1\le m\le N}\chi(y^m).
\lb{3.55}\ee

The following remark shows that $\hat\chi(y)$ is well-defined and is a
rational number.

{\bf Remark 3.16.} By (\ref{3.54}), we have
\be \chi(y^m)
= \sum_{q=0}^{\infty}(-1)^q\dim C_{S^1,\; q}(\Psi_a, \;S^1\cdot u^m)
= \sum_{l=0}^{2n-2}(-1)^{i(y^m)+l}k_l(y^m).\lb{3.56}\ee
Here the first equality follows from Definition 3.1.
The second equality follows from Proposition 3.10 and
Definition 3.11. Hence by (\ref{3.50}) and Proposition 3.13 we have
\bea \hat\chi(y)
&=&\lim_{N\rightarrow\infty}\frac{1}{N}
  \sum_{1\le m\le N\atop 0\le l\le 2n-2}(-1)^{i(y^m)+l}k_l(y^m)\nn\\
&=&\lim_{s\rightarrow\infty}\frac{1}{sK(y)}
  \sum_{1\le m\le K(y),\; 0\le l\le 2n-2\atop 0\le p< s}
  (-1)^{i(y^{pK(y)+m})+l}k_l(y^{pK(y)+m})\nn\\
&=&\frac{1}{K(y)}
  \sum_{1\le m\le K(y)\atop 0\le l\le 2n-2}
  (-1)^{i(y^{m})+l}k_l(y^{m}).\lb{3.57}\eea
Therefore $\hat\chi(y)$ is well defined and is a rational number.
In particular, if all $y^m$s are non-degenerate, then $\nu(y^m)=1$
for all $m\in\N$. Hence the proof of Proposition 3.13 yields $K(y)=2$.
By (\ref{3.46}) we have
\bea k_l(y^m)
    = \left\{\matrix{1, & {\rm if\;\;} i(y^m)-i(y)\in 2\Z \quad {\rm and} \quad l=0  \cr
           0, & {\rm otherwise}. \cr}\right.  \nn\eea
Hence (\ref{3.57}) implies
\be \hat\chi(y)
    = \left\{\matrix{(-1)^{i(y)}, & {\rm if\;\;} i(y^2)-i(y)\in 2\Z,  \cr
           \frac{(-1)^{i(y)}}{2}, & {\rm otherwise}. \cr}\right.  \lb{3.58}\ee

{\bf Remark 3.17.} Note that $k_l(y^m)=0$ for $l\notin [0, \nu(y^m)-1]$
and it can take only values $0$ or $1$ when $l=0$ or $l=\nu(y^m)-1$.
Moreover, the following facts are useful (cf. Lemma 3.10 of \cite{BaL1},
\cite{Cha1} and \cite{MaW1}):

(i) $k_0(y^m)=1$ implies $k_l(y^m)=0$ for $1\le l\le \nu(y^m)-1$.

(ii) $k_{\nu(y^m)-1}(y^m)=1$ implies $k_l(y^m)=0$ for $0\le l\le \nu(y^m)-2$.

(iii) $k_l(y^m)\ge 1$ for some $1\le l\le \nu(y^m)-2$ implies
$k_0(y^m)=k_{\nu(y^m)-1}(y^m)=0$.

(iv) In particular, only one of the $k_l(y^m)$s for $0\le l\le \nu(y^m)-1$ can
be non-zero when $\nu(y^m)\le 3$.

\setcounter{equation}{0}%\setcounter{figure}{0}
\section{Homological vanishing near the origin}%{Section 4}

In Section 3, we have studied nonzero critical points of $\Psi_a$. This section
is devoted to the study of the contribution of the origin to the Morse series of
the functional $\Psi_a$ on $L_0^2(S^1,\R^{2n})$. The main result in this section
is motivated by Theorem 7.1 of \cite{Vit1}, but our proof is different from that
in \cite{Vit1}.

We consider first the distribution of critical values of $\Psi_a$. Note that a
critical point $u_a=u_a(\tau,y)$ of $\Psi_a$ corresponds to a closed
characteristic $(\tau,y)$ on $\Sigma$ by Propositions 2.4 and 2.6. Therefore
the critical value $\Psi_a(u_a)=\Psi_a(u_a(\tau,y))$ is a function of
$(\tau, y)$, $a$ and $\varphi_a$ in Proposition 2.4.

{\bf Proposition 4.1.} {\it Let $H_a(x)=a\varphi(j(x))$ for
$a\in [a_1, a_2]$ be a family of functions given by Proposition 2.4
with the same $\varphi$ satisfying (i)-(iv) of Proposition 2.2.
For $a\in[a_1, a_2]$, suppose $u_a=u_a(\tau,y)\neq 0$ is a critical point of
$\Psi_a$. Then $\Psi_a(u_a)=\Psi_a(u_a(\tau,y))$ is a function of the
period $\tau$, $a$ and $\varphi$. Thus we denote it simply by $\Psi_{a,\tau}$.
Here we ignore the dependence of $\Psi_{a,\tau}$ on $\varphi$. Then we have
the following properties of $\Psi_{a,\tau}$.

(i) $\Psi_{a,\tau}<0$ and $\Psi_{a,\tau}$ is an increasing function of $\tau$
when $a$ is fixed.

(ii) $\Psi_{a,\tau}$ is a strictly decreasing continuous function of $a$ when
$\tau$ is fixed.

(iii) If $a\in per(\Sg)$, then we have $\lim_{\lambda\to a^+}\Psi_{\lambda,a}=0$. }

{\bf Proof.} By Propositions 2.4 and 2.9, we have $u_a=\dot x_a$,
$x_a=\rho_ay(\tau t)$, and then (\ref{2.16}) and (\ref{2.17})) become
\bea
&& \frac{\varphi^\prime(\rho_a)}{\rho_a} = \frac{\tau}{a}, \lb{4.1}\\
&& \Psi_a(u_a) = \frac{1}{2}a\varphi^\prime(\rho_a)\rho_a-a\varphi(\rho_a).
   \lb{4.2}\eea
Note that $\Psi_a(u_a)<0$ follows from Proposition 2.9.
Note that by (\ref{4.1}), $\rho_a$
depends only on the period $\tau$ of the closed characteristic,
$a$ and $\varphi$, hence so does $\Psi_a(u_a)$ by (\ref{4.2}).

For (i), fix an $a\in[a_1, a_2]$.  Note that
by (iii) of Proposition 2.2 and (\ref{4.1}), $\rho_a$ is a decreasing
function of $\tau$. Now let $f(t)=\frac{1}{2}a\varphi^\prime(t)t-a\varphi(t)$.
As in the proof of Proposition 2.9, we obtain $f(0)=0$ and $f^\prime(t)<0$
by (iii) of Proposition 2.2. Hence (i) holds.

For (ii), given $\tau\in per(\Sigma)$ with $\tau<a$,
we may choose a fixed closed characteristic $(\tau,\, y)$ on $\Sigma$.
Differentiating the equation (\ref{4.1}) with respect to $a$ yields
$$ (\varphi^\prime(\rho_a)+a\varphi^{\prime\prime}(\rho_a)\rho^\prime_a)\rho_a
-a\varphi^\prime(\rho_a)\rho^\prime_a=0. $$
Hence we have
\bea \frac{d}{da}\Psi_a(u_a)&=&
\frac{1}{2}\left(\frac{}{}\varphi^\prime(\rho_a)\rho_a
-2\varphi(\rho_a)+a\varphi^{\prime\prime}(\rho_a)\rho_a\rho^\prime_a
-a\varphi^\prime(\rho_a)\rho^\prime_a\right)\nn\\
&=& -\varphi(\rho_a)<0.\nn
\eea
Therefore (ii) holds. Note that in this proof we used the property that $\vf$ is
independent of the choice of $a\in [a_1,a_2]$.

For (iii), we have  $\rho_\lambda\rightarrow 0$ as $\lambda\rightarrow a^+$
by (\ref{4.1}), (i) and (iii) of Proposition 2.2. Hence (iii) holds
by (\ref{4.2}) and (i) of Proposition 2.2. \hfill\hb

The main result in this section is

{\bf Theorem 4.2.} {\it Fix an $a>0$ such that $per(\Sg)\cap (0,a)\not=\emptyset$.
Then there exists an $\vep_0>0$ small such that for any $\vep\in (0,\vep_0]$ we
have
\be H_{S^1,\; q}(\Lambda_a^{\vep}, \;\Lambda_a^{-\vep}) = 0,
 \quad \forall q\le I_0, \lb{4.3}\ee
if $I_0$ is the greatest integer in $\N_0$ such that $I_0<i(\tau, y)$ for all
closed characteristics $(\tau,\, y)$ on $\Sigma$ with $\tau\ge a$.}

{\bf Proof.} Let
$$ \tau(a)=\max \{\tau<a \;|\;\tau\in per(\Sigma)\},
  \qquad \varepsilon_0=-\frac{1}{2}\Psi_{a,\tau(a)}.$$
Then by (i) of Proposition 4.1, there are no critical values of $\Psi_a$ in
the interval $[-\varepsilon_0,\; \varepsilon_0]$ except $0$. Hence we have
\be H_{S^1,\; q}(\Lambda_a^{\vep}, \;\Lambda_a^{-\vep})
  \cong H_{S^1,\; q}(\Lambda_a^{\vep_0}, \;\Lambda_a^{-\vep_0}),
  \quad \forall q\in\Z,\; \vep\in(0,\,\vep_0]. \lb{4.4}\ee
In the following we assume $\vep\in(0,\,\vep_0]$.

Note that by the same proof of Proposition 3.2,
$H_{S^1,\; q}(\Lambda_a^\varepsilon, \;\Lambda_a^{-\varepsilon})$ is
independent of the choice of $\varphi_a$ in $H_a(x)=a\varphi_a(j(x))$
which satisfies (i)-(iv) of Proposition 2.2. Hence we can choose
$\varphi_a\equiv\varphi$ for any function $\vf$ satisfying (i)-(iv)
of Proposition 2.2.

The rest part of the proof of this theorem is carried out in three steps.

{\bf Step 1.} Claim: For every $b>a$ there exists $\hat{\vep}_b\in(0,\;\vep_0]$
such that
\be H_{S^1,\; q}(\Lambda_a^{\vep}, \;\Lambda_a^{-\vep})
  \cong H_{S^1,\; q}(\Lambda_b^{\vep}, \;\Lambda_b^{-\vep}),
  \quad\forall q\le I_0 \;\,{\rm and}\;\, \vep\in (0,\hat{\vep}_b]. \lb{4.5}\ee

In fact, by the above second paragraph, we may choose $\varphi$ such
that $H_c(x)=c\varphi(j(x))$ satisfies Proposition 2.4 for all
$c\in[a, b]$ with a fixed $\varphi$.

By Lemma 3.4, we can choose a family of finite dimensional
approximations $h_c: X\rightarrow X^\perp$ and consider
the functions $\psi_c(g)=\Psi_c(g+h_c(g))$ on the finite dimensional
manifold $X$. Moreover, we have
$H_{S^1,\; q}(\Lambda_c^\varepsilon, \;\Lambda_c^{-\varepsilon})
\cong H_{S^1,\; q}(\widetilde{\Lambda}_c^\varepsilon, \;
\widetilde{\Lambda}_c^{-\varepsilon})$ for any $\varepsilon>0$
by (iv) of Lemma 3.4. Hence in order to prove (\ref{4.5}),
it suffices to prove
\be H_{S^1,\; q}(\widetilde{\Lambda}_a^\varepsilon, \;
\widetilde{\Lambda}_a^{-\varepsilon})
\cong H_{S^1,\; q}(\widetilde{\Lambda}_b^\varepsilon, \;
\widetilde{\Lambda}_b^{-\varepsilon}),\lb{4.6}\ee
for $\varepsilon\in (0, \varepsilon_0)$ sufficiently small.
Clearly, it suffices to prove that for any $c\in[a, b]$
there exists $\delta,\,\varepsilon^\prime>0$ such that
\be H_{S^1,\; q}(\widetilde{\Lambda}_{c_1}^\varepsilon, \;
\widetilde{\Lambda}_{c_1}^{-\varepsilon})
\cong H_{S^1,\; q}(\widetilde{\Lambda}_{c_2}^\varepsilon, \;
\widetilde{\Lambda}_{c_2}^{-\varepsilon}),\quad\forall q\le I_0.\lb{4.7}\ee
for any $c_1, c_2\in[c-\delta, c+\delta]$ and $\varepsilon\in(0,\,\varepsilon^\prime]$.
In the following, we fix a $c\in[a, b]$.

We have two cases.

{\bf Case 1.} $c\notin per(\Sigma)$.

In this case, since $per(\Sigma)$ is a discrete subset of $\R^+$ by
Definition 2.1 and the assumption (F), we can find $\delta>0$ such that
$[c-\delta, c+\delta]\cap per(\Sigma)=\emptyset$. Hence nonzero critical
points of $\psi_\lambda$ are precisely those closed characteristics
$(\tau, y)$ on $\Sigma$ with period $\tau<c$ for all
$\lambda\in[c-\delta, c+\delta]$. Let
$$ \tau_0=\max \{\tau<c-\delta\;|\;\tau\in per(\Sigma)\}. $$
Then by Proposition 4.1 and the definition of $\tau_0$, we have
$$ \psi_\lambda(g_\lambda)\le \Psi_{\lambda, \tau_0}\le\Psi_{c-\delta, \tau_0}<0 $$
for all nonzero critical points $g_\lambda$ of $\psi_\lambda$.
Let $\varepsilon^\prime=-\frac{1}{2}\Psi_{c-\delta, \tau_0}$.
Then $\pm \varepsilon$ for $\varepsilon\in (0,\,\varepsilon^\prime]$
are regular values of $\psi_\lambda$ for all $\lambda\in[c-\delta, c+\delta]$.
Moreover, by (\ref{2.15}) and the definition of $\psi_{\lm}$ in Lemma 3.4,
we have $\psi_\lambda(g)$ goes
to $+\infty$ as $\|g\|$ goes to $+\infty$. Hence we can choose $R>0$ to be
sufficiently large such that $\widetilde{\Lambda}_{\lambda}^{\vep}\subset B_R(0)$
for all $\lambda\in[c-\delta, c+\delta]$, where $B_R(0)$ is the open ball in $X$
centered at the origin with the radius $R$. Then we have
$\psi_\lambda(g)$ is continuous with respect to $\lambda$ uniformly for
$g\in\bigcup_{\lambda\in[c-\delta, c+\delta]}\widetilde{\Lambda}_\lambda^\varepsilon$.
Clearly, it follows from Lemma 3.4 that both $\psi_\lambda(g)$ and
$\psi_\lambda^\prime(g)$ are continuous for
$(\lambda,\, g)\in [c-\delta, c+\delta]\times X$.
Hence we can apply a slightly stronger version of Exercise 8.4 on p.203 of
\cite{MaW1} or Theorem 1.5.6 on p.53 of \cite{Cha1} to obtain (\ref{4.7}).

{\bf Case 2.} $c\in per(\Sigma)$.

In this case, let
$$ \tau_0=\max \{\tau<c \;|\;\tau\in per(\Sigma)\}, \qquad
   \tau_1=\min \{\tau>c \;|\;\tau\in per(\Sigma)\}.  $$
Let $\delta<\frac{1}{2}\min\{c-\tau_0, \;\tau_1-c\}$ to be
determined later. Then nonzero critical points of $\psi_\lambda$
for $\lambda\in[c-\delta, c]$ consists of closed characteristics
$(\tau, y)$ on $\Sigma$ with period $\tau<c$, and nonzero critical
points of $\psi_{\lm}$ for $\lm\in (c, c+\delta]$ consists of
closed characteristics $(\tau,y)$ on $\Sigma$ with period $\tau\le
c$. Then by Proposition 4.1, we have
$$ \psi_{\lm}(g_\lm)\le \Psi_{\lm,\,\tau_0}\le \Psi_{c-\delta,\,\tau_0}<0  $$
for all nonzero critical points $g_{\lm}$ of $\psi_{\lm}$
with period $\tau<c$ when $\lm\in [c-\delta, c+\delta]$.
Note that by (ii) and (iii) of Proposition 4.1, the
critical value $\Psi_{\lm,\,c}$ of $\psi_{\lm}$ for $\lm\in (c, c+\delta]$
with period $\tau=c$ is close to $0$ when $\delta$ is small. Specially by
Proposition 4.1, from the fact
$\Psi_{c,\,\tau_0}<\frac{1}{2}\Psi_{c,\,\tau_0}<0
    =\lim_{\lambda\rightarrow c^+}\Psi_{\lambda,\,c}$
we can choose $\delta>0$ so small such that the following relation holds:
$$ \frac{1}{2}\Psi_{c,\,\tau_0} < \frac{1}{2}\Psi_{c-\delta, \,\tau_0}
    < \Psi_{c+\delta,\,c} < 0=\lim_{\lambda\rightarrow c^+}\Psi_{\lambda,\,c}. $$
Let $\vep'=-\frac{1}{2}\Psi_{c-\delta, \tau_0}$. Therefore by our above
discussion $\pm \vep'$ are regular values of $\psi_{\lm}$ for all
$\lm\in[c-\delta, c+\delta]$.

\begin{figure}
\begin{center}
\resizebox{9cm}{6.9cm}{\includegraphics*[0cm,0cm][18.06cm,13.55cm]{WHL-fig.eps}}
%%{Figure 4.1}
\caption{Near the critical value $0$}
\end{center}
\end{figure}

More precisely, we have Figure 4.1, where the $\mu$-axis denote the
critical values of $\psi_{\lm}$, $f_0(\lm)=\Psi_{\lm,\tau_0}$ for
$\lm\in[c-\delta, c+\delta]$ and $f_1(\lm)=\Psi_{\lm,c}$ for
$\lm\in(c, c+\delta]$. Both $f_0$ and $f_1$ are decreasing functions
in $\lm$ by (ii) of Proposition 4.1, and $\lim_{\lm\to c^+}f_1(\lm)=0$
by (iii) of Proposition 4.1. Moreover, by (i) of Proposition 4.1, the
interior of the shaded part in the Figure 4.1 contains no critical values
of $\psi_{\lm}$ for $\lm\in [c-\delta, c+\delta]$. Therefore $-\vep'$
is a common regular value of those $\psi_{\lm}$s. Since all critical
values of $\psi_{\lm}$ are non-positive, $\vep'$ is a common regular
value for all of $\psi_{\lm}$s too.

Hence by the same proof of Case 1, we have
\be H_{S^1,\;q}(\widetilde{\Lambda}_{c_1}^{\varepsilon^\prime}, \;
 \widetilde{\Lambda}_{c_1}^{-\varepsilon^\prime})
 \cong H_{S^1,\; q}(\widetilde{\Lambda}_{c_2}^{\varepsilon^\prime}, \;
 \widetilde{\Lambda}_{c_2}^{-\varepsilon^\prime}),\quad \forall q\in \Z
 \lb{4.8}\ee
for any $c_1, c_2\in[c-\delta, c+\delta]$.

Since $0$ is the only critical value of $\psi_{\lm}$ in $[-\vep', \vep']$
for $\lm\in [c-\delta, c]$, for any $\vep\le\vep'$ we have
\be H_{S^1, q}(\widetilde{\Lm}_{\lm}^{\vep'}, \widetilde{\Lm}_{\lm}^{-\vep'})
 \cong H_{S^1, q}(\widetilde{\Lm}_{\lm}^{\vep}, \widetilde{\Lm}_{\lm}^{-\vep}),
 \quad \forall q\in \Z. \lb{4.9}\ee

But critical values of $\psi_{\lm}$ in $[-\vep', \vep']$ for $\lm\in(c, c+\dl]$
are precisely $0$ and $\Psi_{\lm,\,c}$ as indicated in the Figure 4.1.

If $\varepsilon\in(-\Psi_{\lm,\,c},\, \vep']$, then the interval
$[-\vep',\,-\vep]$ contains no critical values of $\Psi_\lm$. Hence
(\ref{4.9}) remains true for these $\vep$ and $\lm$.

If $\vep\in(0,\,-\Psi_{\lm,\,c}]$, we consider the exact sequence of the triple
$(\widetilde{\Lm}_{\lm}^{\vep}, \widetilde{\Lm}_{\lm}^{-\vep},
   \widetilde{\Lm}_{\lm}^{-\vep'})$:
\be
H_{S^1, q}(\widetilde{\Lm}_{\lm}^{-\vep}, \widetilde{\Lm}_{\lm}^{-\vep'})
 \rightarrow H_{S^1, q}(\widetilde{\Lm}_{\lm}^{\vep}, \widetilde{\Lm}_{\lm}^{-\vep'})
 \rightarrow H_{S^1, q}(\widetilde{\Lm}_{\lm}^{\vep}, \widetilde{\Lm}_{\lm}^{-\vep})
 \rightarrow H_{S^1, q-1}(\widetilde{\Lm}_{\lm}^{-\vep}, \widetilde{\Lm}_{\lm}^{-\vep'}).
          \lb{4.10} \ee
Since $\Psi_{\lambda, c}$ is the unique critical
value of $\psi_\lambda$ in $[-\varepsilon^\prime, -\varepsilon]$,
as in Lemma 1.4.2 of \cite{Cha1}, we have
\bea
H_{S^1,\; \ast}(\widetilde{\Lambda}_{\lambda}^{-\varepsilon},\;
\widetilde{\Lambda}_{\lambda}^{-\varepsilon^\prime })
&\cong&
\bigoplus_{i=1}^{l}C_{S^1,\;\ast}(\psi_\lambda, \;S^1\cdot g_{\lambda,\, i})\nn\\
&\cong&
\bigoplus_{i=1}^{l}C_{S^1,\;\ast}(\Psi_\lambda,\; S^1\cdot u_{\lambda, \,i}).\lb{4.11}
\eea
Here $g_{\lambda, \,i}$ and $u_{\lambda,\, i}$ denote the critical
points of $\psi_\lambda$ and $\Psi_\lambda$ with critical value
$\Psi_{\lambda,\, c}$ respectively. The second isomorphism follows
from (iv) of Lemma 3.4. By Proposition 3.10 and the definition of $I_0$, we have
$$ C_{S^1,\;q}(\Psi_\lambda, \;S^1\cdot u_{\lambda,\, i})
\cong 0,\quad \forall q\le I_0,\; 1\le i\le l.$$
Hence (\ref{4.11}) yields
\be H_{S^1,\; q}(\widetilde{\Lambda}_{\lambda}^{-\varepsilon},\;
             \widetilde{\Lambda}_{\lambda}^{-\varepsilon^\prime })
\cong 0,\quad \forall q\le I_0.\lb{4.12}\ee
Then by (\ref{4.10}) we have
\be H_{S^1,\; q}(\widetilde{\Lambda}_{\lambda}^{\varepsilon},\;
             \widetilde{\Lambda}_{\lambda}^{-\varepsilon^\prime })
\cong
H_{S^1,\; q}(\widetilde{\Lambda}_{\lambda}^{\varepsilon},\;
             \widetilde{\Lambda}_{\lambda}^{-\varepsilon }),
\quad\forall q\le I_0.\lb{4.13}\ee
Since $\psi_\lambda$ has no critical value
in $[\varepsilon, \varepsilon^\prime]$, we have
\be H_{S^1,\; q}(\widetilde{\Lambda}_{\lambda}^{\varepsilon},\;
             \widetilde{\Lambda}_{\lambda}^{-\varepsilon^\prime })
\cong
H_{S^1,\; q}(\widetilde{\Lambda}_{\lambda}^{\varepsilon^\prime},\;
             \widetilde{\Lambda}_{\lambda}^{-\varepsilon^\prime }),
\quad\forall q\in\Z. \lb{4.14}\ee
Combining (\ref{4.8}), (\ref{4.9}), (\ref{4.13}) and (\ref{4.14})
we obtain (\ref{4.7}). The proof of Step 1 is complete.

{\bf Step 2.} Claim:
\be H_{S^1,\; q}(\Lambda_b^\varepsilon, \;\Lambda_b^{-\varepsilon})
\cong 0,\quad\forall q\le I_0  \lb{4.15}\ee
holds for some $b>a$ large enough and some $\varepsilon\in (0, \hat{\vep}_b]$
sufficiently small.

In fact, the proof is a modification of that of Theorem 3.8 in
\cite{Eke1} to the $S^1$-equivariant case. Considering $\psi_b$, we
assume $b\notin per(\Sigma)$ and will determine $b$ later.

Firstly we approximate $\psi_b$ by an $S^1$-invariant $C^2$
function $\widehat{\psi}$ satisfying the following conditions:

(i) $\widehat{\psi}$ has the same critical points as $\psi_b$ outside
a neighborhood $\Omega$ of $0$. Hence $\widehat{\psi}$ contains all
nonzero critical points of $\psi_b$ as its critical points.

(ii) Each critical orbit $S^1\cdot g$ of $\widehat{\psi}$ contained
in $\Omega$ is non-degenerate and $\widehat{\psi}$ has Morse index
$m^-(g)>I_0$ at the critical point $g$.

More precisely, we construct $\widehat{\psi}$ as follows.

Following p.46 of \cite{Eke1} and Proposition 2.5, we can approximate
$G_b$ by  a $C^2$ strictly convex function $\widetilde{G}$ such that
\bea \widetilde{G}(x)&=&G_b(x),\quad {\rm for}\quad |x|\ge \varrho_1, \lb{4.16}\\
(\widetilde{G}^{\prime\prime}(x)\xi, \xi)
&\le& \frac{1}{br_0}|\xi|^2,
\quad {\rm for}\quad |x|\le \varrho_2, \; \forall \xi\in\R^{2n},
\lb{4.17}   \eea
where $\varrho_2>\varrho_1>0$ can be chosen as small as we want
and $r_0$ is given by Proposition 2.5. Now we define a
functional $\widetilde{\Psi}$ on $L_0^2(S^1, \R^{2n})$
\be \widetilde{\Psi}(u)=\int_0^1\left(\frac{1}{2}Ju\cdot Mu+
\widetilde{G}(-Ju)\right)dt. \lb{4.18}\ee
Denote by $\widetilde{H}$ the Fenchel transform of $\widetilde{G}$.
Then critical points of $\widetilde{\Psi}$ correspond to $1$-periodic
solutions of the equation $\dot x=J\widetilde{H}^\prime(x)$.
Moreover, by choosing $\varrho_1$ small enough,
nonzero critical points of $\Psi_b$ are also critical points of
$\widetilde{\Psi}$. Other critical points $u$ of $\widetilde{\Psi}$
must satisfy $\|u\|_{C^0(S^1, \R^{2n})}<\varrho_2$. Hence such a
critical point $u$ has index
\be i(u)\ge 2n\left[\frac{br_0}{2\pi}\right]\equiv I_1(b),
\lb{4.19}\ee
by Definition 3.3, (\ref{4.17}) and Proposition 1.4.14 on p.32 of
\cite{Eke3}. Now we can fix the $b$ in the Claim (\ref{4.15}) to
satisfy $I_1(b)>I_0$.

For any $\kappa\in\R$, we denote by
\be \Theta^\kappa=\{u\in L_0^2(S^1, \R^{2n}) \;|\;
\widetilde{\Psi}(u)\le\kappa\}. \lb{4.20}\ee
By choosing $\|\widetilde{\Psi}-\Psi_b\|_{C^1(L_0^2(S^1, \R^{2n}),\, \R)}$
to be small enough, we can fix an $\vep\in (0,\hat{\vep}_b)$ such that
$(-2\vep, 0)$ contains no critical value of $\Psi_b$ and $\pm\vep$
are regular values of $\widetilde{\Psi}$.
Then using a slightly stronger version of Exercise 8.4 of \cite{MaW1}
or Theorem 1.5.6 on p.53 of \cite{Cha1} with continuous dependence on
the parameter, we obtain
\be H_{S^1,\; q}(\Lm_b^{\vep},\; \Lm_b^{-\vep})
 \cong H_{S^1,\; q}(\Theta^{\vep}, \;\Theta^{-\vep}),
 \quad\forall q\in\Z.  \lb{4.21}\ee
By Lemma 3.4, we can choose a finite dimensional approximation
$\widetilde{h}: X\rightarrow X^\perp$ and consider the function
$\widetilde{\psi}(g)=\widetilde{\Psi}(g+\widetilde{h}(g))$ on the
finite dimensional manifold $X$. We have
\be H_{S^1,\,q}(\Theta^{\vep}, \;\Theta^{-\vep})
  \cong H_{S^1,\,q}(\widetilde{\Theta}^{\vep}, \;\widetilde{\Theta}^{-\vep}),
  \quad\forall q\in\Z. \lb{4.22}\ee
by (iv) of Lemma 3.4, where
$\widetilde{\Theta}^\kappa=\{g\in X\;|\; \widetilde{\psi}(g)\le \kappa\}$
for $\kappa\in\R$ and any critical  point $g$ of $\widetilde{\psi}$
with critical value in $[-\varepsilon, \varepsilon]$ has Morse index
$m^-(g)\ge I_1(b)$ by (iii) of Lemma 3.4 and (\ref{4.19}). By (iii)
of Lemma 3.4, we have $\widetilde{\psi}\in C^2(X, \;\R)$ is $S^1$-invariant,
and the $S^1$-action is $C^\infty$ on $X$. Hence by  the
Density Lemma of \cite{Was1}, $\widetilde{\psi}$ can be $C^2$ approximated
by a smooth $S^1$-invariant function $\widehat{\psi}$ whose critical
orbits $S^1\cdot g$ are non-degenerate when $\widehat{\psi}(g)\in [-\vep,\vep]$,
i.e., $\widehat{\psi}$ is a Morse function there, and any critical point $g$ of
$\widehat{\psi}$ with critical values in $[-\varepsilon, \varepsilon]$ has Morse
index $m^-(g)\ge I_1(b)$. This finish the construction of $\widehat{\psi}$.

When $\widehat{\psi}$ is sufficiently close to $\widetilde{\psi}$, we have
by a slightly stronger version of Exercise 8.4 of \cite{MaW1} or Theorem
1.5.6 on p.53 of \cite{Cha1} again,
\be H_{S^1,\; q}(\widetilde{\Theta}^\varepsilon, \;\widetilde{\Theta}^{-\varepsilon})
\cong H_{S^1,\; q}(\Delta^\varepsilon, \;\Delta^{-\varepsilon}),
\quad\forall q\in\Z, \lb{4.23}\ee
where $\Delta^\kappa=\{g\in X\;|\; \widehat{\psi}(g)\le \kappa\}$ for $\kappa\in\R$.
Now by the Thom isomorphism (cf. p.77 of \cite{Cha1}), we have
\be C_{S^1,\; q}(\widehat{\psi},\; S^1\cdot g)\cong
H_{q-m^-(g)}(BG_g,\; \theta),\lb{4.24}\ee
where $G_g$ is the isotropy group of the critical orbit $S^1\cdot g$
and $\theta$ is the orientation bundle of the negative bundle
of $\widehat{\psi}^{\prime\prime}(g)$. Hence we have
\be C_{S^1,\; q}(\widehat{\psi},\; S^1\cdot g)\cong 0,\quad
  \forall q<I_1(b), \lb{4.25}\ee
for any critical point $g$ of $\widehat{\psi}$
with critical values in $[-\varepsilon, \varepsilon]$.
Hence by the Morse inequality, we have
\be \sum_{i=1}^{l}\dim C_{S^1,\; q}(\widehat{\psi},\; S^1\cdot g_i)
\ge \dim H_{S^1,\; q}(\Delta^\varepsilon, \;\Delta^{-\varepsilon}),
\quad\forall q\in\Z, \lb{4.26}\ee
where we denote the critical orbits of $\widehat{\psi}$
with critical values in $[-\varepsilon, \varepsilon]$
by $\{S^1\cdot g_1, \ldots, S^1\cdot g_l\}$.
Now combining (\ref{4.21})-(\ref{4.23}), (\ref{4.25})
and (\ref{4.26}), we obtain the claim (\ref{4.15}).

{\bf Step 3.} Now (\ref{4.5}) of Step 1 and (\ref{4.15}) of Step 2 yield an
$\vep\in (0,\hat{\vep}_b]$ for some $b>a$ large enough such that (\ref{4.3})
holds for this $\vep$. Then by (\ref{4.4}) we obtain (\ref{4.3}) for all
$\vep\in (0,\ep_0]$ and then the proof of Theorem 4.2 is complete.
\hfill\hb

\setcounter{equation}{0}%\setcounter{figure}{0}
\section{ Proof of the Theorem 1.2}%{Section 5}

In this section, we give a proof for the Theorem 1.2 with
$H_a(x)=a\vf_a(j(x))$, where $\vf_a$ satisfies (i)-(iv) of Proposition
2.2.

Let $\Psi_a$ be a functional defined by (\ref{2.11}) for
some $a\in\R$ large enough and let $\varepsilon>0$ be small
enough such that $[-\varepsilon, 0)$ contains no critical
values of $\Psi_a$. We consider the exact sequence of the space pair
$(\Lm_a^{\infty},\;\Lm_a^{-\vep})$:
\be
H_{S^1,\; q+1}(\Lm_a^{\infty},\;\Lm_a^{-\vep} )
 \rightarrow H_{S^1,\; q}(\Lm_a^{-\vep})
 \rightarrow H_{S^1,\; q}(\Lm_a^{\infty})
 \rightarrow H_{S^1,\; q}(\Lm_a^{\infty},\;\Lm_a^{-\vep})
\lb{5.1}\ee
for any $q\in\Z$. Let $I_0\in\N_0$ be given by Theorem 4.2.
Note that by Proposition 4.1, there are no critical values of
$\Psi_a$ in $(0, +\infty)$. Hence by Theorem 4.2 we have
\be
H_{S^1,\; q}(\Lm_a^\infty,\;\Lm_a^{-\vep})
\cong H_{S^1,\; q}(\Lm_a^{\vep},\;\Lm_a^{-\vep})
\cong 0, \quad \forall q\le I_0.\lb{5.2}\ee
Therefore (\ref{5.1}) implies
\bea
H_{S^1,\; q}(\Lambda_a^{-\varepsilon} )
\cong
H_{S^1,\; q}( \Lambda_a^\infty)
\cong
H_q(CP^\infty), \quad \forall q<I_0.\lb{5.3}\eea
The second isomorphism follows since
$\Lambda_a^\infty=L^2_0(S^1,\;\R^{2n})$
is $S^1$-equivariantly homotopic to a single point.

Let $X$ be an $S^1$-space such that the Betti numbers
$b_i(X)=\dim H_{S^1,\;i}(X;\;\Q)$ are finite for all $i\in \Z$.
As usual the $S^1$-equivariant {\it Poincar\'e series} of $X$ is
defined by the formal power series $P(X)(t)=\sum_{i=0}^{\infty}b_i(X)t^i$.
Note that by Proposition 2.7, $\Psi_a$ is bounded from below on
$L^2_0(S^1,\;\R^{2n})$. Hence the $S^1$-equivariant {\it Morse
series} $M(t)$ of the functional $\Psi_a$ on the space
$\Lambda_a^{-\varepsilon}$ is defined as usual by
\bea
M(t)=\sum_{q\ge 0,\;1\le j\le p} \dim C_{S^1,\;q}(\Psi_a, \;S^1\cdot
v_j)t^q,\lb{5.4} \eea
where we denote by $\{S^1\cdot v_1, \ldots,
S^1\cdot v_p\}$ the critical orbits of $\Psi_a$ with critical values
less than $-\varepsilon$.
Then the Morse inequality in the equivariant sense yields a formal power
series $Q(t)=\sum_{i=0}^\infty q_it^i$ with nonnegative integer
coefficients $q_i$ such that
\bea  M(t)=P(t)+(1+t)Q(t), \lb{5.5}\eea
where $P(t)\equiv P(\Lambda_a^{-\varepsilon})(t)$.
For a formal power series
$R(t)=\sum_{i=0}^\infty r_it^i$, we denote by $R^L(t)=\sum_{i=0}^L
r_i t^i$ for $L\in\N$ the corresponding truncated polynomial. Using
this notation, (\ref{5.5}) becomes
\be (-1)^Lq_L=M^L(-1)-P^L(-1), \quad \forall L\in\N.\lb{5.6}\ee

Now we can give the following

{\bf Proof of Theorem 1.2.} Firstly we choose $\Psi_a$ as above and
denote by $\{u_1, \ldots, u_k\}$ the  critical points of $\Psi_a$
corresponding to $\{y_1,\ldots, y_k\}$. Note that $v_1,\ldots,v_p$ in
(\ref{5.4}) are iterations of $u_1,\ldots,u_k$. Since
$C_{S^1,\;q}(\Psi_a, \;S^1\cdot u_j^m)$ can be non-zero
only for $q=i(u_j^m)+l$ with $0\le l\le 2n-2$ by Propositions
3.5 and 3.10, the formal Poincar\'e series (\ref{5.4}) becomes
\be M(t)=\sum_{1\le j\le k,\; 0\le l\le 2n-2 \atop 1\le m_j<a/\tau_j}
              k_l(u_j^{m_j})t^{i(u_j^{m_j})+l}
  = \sum_{1\le j\le k,\; 0\le l\le 2n-2 \atop 1 \le m_j\le K_j,\; sK_j+m_j<a/\tau_j}
               k_l(u_j^{m_j})t^{i(u_j^{sK_j+m_j})+l}, \lb{5.7}\ee
where $K_j=K(u_j)$ and $s\in\N_0$. The last equality follows from Proposition 3.13.
Let $I=I_0-2$, where $I_0$ is given by (\ref{5.3}) and consider the
truncated polynomials $M^{I}(t)$ and $P^{I}(t)$.

Write $M(t)=\sum_{h=0}^{\infty}w_ht^h$ and
$P^{I}(t)=\sum_{h=0}^{I}b_ht^h$. Then we have
\be w_h\ = \sum_{1\le j\le k,\; 0\le l\le 2n-2 \atop 1 \le m\le K_j}
              k_l(u_j^m)\,^\#\{s\in\N_0\,|\,i(u_j^{sK_j+m})+l=h\},
\quad \forall h\le I+1. \lb{5.8}\ee
Note that the right hand side of (\ref{5.7}) contains only those terms
satisfying $sK_j+m_j<\frac{a}{\tau_j}$. Thus (\ref{5.8}) holds only for
$h\le I+1$ by (\ref{5.7}).

{\bf Claim 1.} {\it $w_h\le C$ for $h\le I+1$ with $C$ being
independent of $a$}.

In fact, we have
\bea
^\#\{s\in\N_0 &|& i(u_j^{sK_j+m})+l=h \}\nn\\
&=&\;^\#\{s\in\N_0 \;| \; i(u_j^{sK_j+m})+l=h,\;
                          |i(u_j^{sK_j+m})-(sK_j+m)\hat{i}(u_j)|\le 2n\} \nn\\
&\le &\;^\#\{s\in\N_0 \;| \;|h-l-(sK_j+m)\hat{i}(u_j)|\le 2n\}  \nn\\
&=&\;^\#\left\{s\in\N_0 \; \left|\;\frac{}{}\right.
      \;\frac{h-l-2n-m\hat{i}(u_j)}{K_j\hat{i}(u_j)}\le s
       \le \frac{h-l+2n-m\hat{i}(u_j)}{K_j\hat{i}(u_j)}\right\}  \nn\\
&\le&\; \frac{4n}{K_j\hat{i}(u_j)}+2,  \lb{5.9}\eea
where the first equality follows from the fact
\be|i(u_j^m)-m\hat{i}(u_j)|\le 2n,\quad \forall m\in\N,\; 1\le j\le k,\lb{5.10}\ee
which follows from Theorems 10.1.2 and 15.1.1 of \cite{Lon4}. Hence Claim 1 holds.

We estimate next $M^I(-1)$. By (\ref{5.8}) we obtain
\bea M^I(-1)
&=& \sum_{h=0}^I w_h(-1)^h   \nn\\
&=& \sum_{1\le j\le k,\; 0\le l\le 2n-2 \atop 1 \le m\le K_j}
            (-1)^{i(u_j^m)+l}k_l(u_j^m)
              \,^\#\{s\in\N_0 \,|\, i(u_j^{sK_j+m})+l\le I\}.
\lb{5.11}\eea
Here the second equality holds by (\ref{3.50}).

{\bf Claim 2.} {\it There is a real constant $C^\prime>0$ independent of $a$
such that
\be \left|M^I(-1)-\sum_{1\le j\le k,\; 0\le l\le 2n-2 \atop 1 \le m\le K_j}
            (-1)^{i(u_j^m)+l}k_l(u_j^m)\frac{I}{K_j\hat{i}(y_j)}\right|
               \le C^\prime,  \lb{5.12}\ee
where the sum in the left hand side of (\ref{5.12}) equals to
$\;I\sum_{1\le j\le k}\frac{\hat\chi(y_j)}{\hat i(y_j)}\;$ by (\ref{3.57}).}

In fact, we have the estimates
\bea
^\#\{s\in\N_0 &|& i(u_j^{sK_j+m})+l\le I\}   \nn\\
&=&\;^\#\{s\in\N_0 \;| \; i(u_j^{sK_j+m})+l\le I,\;
             |i(u_j^{sK_j+m})-(sK_j+m)\hat{i}(u_j)|\le 2n\}  \nn\\
&\le&\;^\#\{s\in\N_0 \;| \;0\le (sK_j+m)\hat{i}(u_j)\le I-l+2n\}  \nn\\
&=&\;^\#\left\{s\in\N_0 \; \left |\;\frac{}{}\right.
    \;0\le s\le \frac{I-l+2n-m\hat{i}(u_j)}{K_j\hat{i}(u_j)}\right\}  \nn\\
&\le&\; \frac{I-l+2n}{K_j\hat{i}(u_j)}+1.  \nn\eea
On the other hand, we have
\bea
^\#\{s\in\N_0 &|& i(u_j^{sK_j+m})+l\le I\}  \nn\\
&=&\;^\#\{s\in\N_0 \;| \; i(u_j^{sK_j+m})+l\le I,\;
               |i(u_j^{sK_j+m})-(sK_j+m)\hat{i}(u_j)|\le 2n\}  \nn\\
&\ge&\;^\#\{s\in\N_0\;|\;i(u_j^{sK_j+m})\le(sK_j+m)\hat{i}(u_j)+2n\le I-l\} \nn\\
&\ge&\;^\#\left\{s\in\N_0 \; \left |\;\frac{}{}\right.
   \;0\le s\le \frac{I-l-2n-m\hat{i}(u_j)}{K_j\hat{i}(u_j)}\right\}  \nn\\
&\ge&\;\frac{I-l-2n}{K_j\hat{i}(u_j)}-2,  \nn\eea
where $m\le K_j$ is used.
Combining these two estimates together with (\ref{5.11}), we obtain (\ref{5.12}).

Note that all coefficients in (\ref{5.5}) are nonnegative, hence by Claim 1,
we have $q_I\le w_{I+1}\le C$. By (\ref{5.3}), we have
$P^I(t)=\sum_{0\le h\le\frac{I}{2}}t^{2h}$.

By (\ref{5.6}), we have
\be (-1)^Iq_I=M^I(-1)-P^I(-1)=M^I(-1)-\left(\left[\frac{I}{2}\right]+1\right).\lb{5.13} \ee
By Theorem 1.7.7 of \cite{Eke3} or Lemma 15.3.2 of \cite{Lon4}, we have $\hat i(y_j)>2$
for $1\le j\le k$.
Hence $i(m\tau_j,\, y_j)\equiv i(y_j^m)\rightarrow\infty$ as
$m\rightarrow\infty$ for $1\le j\le k$.
Now we let $a\rightarrow+\infty$, then $I=I_0-2\rightarrow+\infty$ in Theorem 4.2.
Note that by Claims 1 and 2, the constants $C$ and $C^\prime$ are independent of
$a$. Hence dividing both sides of (\ref{5.13}) by $I$ and letting $I$ tending
to infinity yield
$$ \lim_{I\to\infty}\frac{1}{I}M^I(-1) = \frac{1}{2}. $$
Hence (\ref{1.3}) holds by (\ref{5.12}).  \hfill\hb

\setcounter{equation}{0}%\setcounter{figure}{0}
\section{ Proofs of the Theorems 1.1 and 1.4}%{Section 6}

In this section, we prove Theorems 1.1 and 1.4 based on Theorem 1.2 and the
index iteration theory developed by Y. Long and his coworkers.

\subsection{A brief review on an index theory for symplectic paths}

In this subsection, we recall briefly an index theory for symplectic paths.
All the details can be found in \cite{Lon4}.

As usual, the symplectic group $\Sp(2n)$ is defined by
$$ \Sp(2n) = \{M\in {\rm GL}(2n,\R)\,|\,M^TJM=J\}, $$
whose topology is induced from that of $\R^{4n^2}$. For $\tau>0$ we are interested
in paths in $\Sp(2n)$:
$$ \P_{\tau}(2n) = \{\ga\in C([0,\tau],\Sp(2n))\,|\,\ga(0)=I_{2n}\}, $$
which is equipped with the topology induced from that of $\Sp(2n)$. The
following real function was introduced in \cite{Lon2}:
$$ D_{\om}(M) = (-1)^{n-1}\ol{\om}^n\det(M-\om I_{2n}), \qquad
          \forall \om\in\U,\, M\in\Sp(2n). $$
Thus for any $\om\in\U$ the following codimension $1$ hypersurface in $\Sp(2n)$ is
defined in \cite{Lon2}:
$$ \Sp(2n)_{\om}^0 = \{M\in\Sp(2n)\,|\, D_{\om}(M)=0\}.  $$
For any $M\in \Sp(2n)_{\om}^0$, we define a co-orientation of $\Sp(2n)_{\om}^0$
at $M$ by the positive direction $\frac{d}{dt}Me^{t\ep J}|_{t=0}$ of
the path $Me^{t\ep J}$ with $0\le t\le 1$ and $\ep>0$ being sufficiently
small. Let
\bea
\Sp(2n)_{\om}^{\ast} &=& \Sp(2n)\bs \Sp(2n)_{\om}^0,   \nn\\
\P_{\tau,\om}^{\ast}(2n) &=&
      \{\ga\in\P_{\tau}(2n)\,|\,\ga(\tau)\in\Sp(2n)_{\om}^{\ast}\}, \nn\\
\P_{\tau,\om}^0(2n) &=& \P_{\tau}(2n)\bs  \P_{\tau,\om}^{\ast}(2n).  \nn\eea
For any two continuous arcs $\xi$ and $\eta:[0,\tau]\to\Sp(2n)$ with
$\xi(\tau)=\eta(0)$, it is defined as usual:
$$ \eta\ast\xi(t) = \left\{\matrix{
            \xi(2t), & \quad {\rm if}\;0\le t\le \tau/2, \cr
            \eta(2t-\tau), & \quad {\rm if}\; \tau/2\le t\le \tau. \cr}\right. $$
Given any two $2m_k\times 2m_k$ matrices of square block form
$M_k=\left(\matrix{A_k&B_k\cr
                                C_k&D_k\cr}\right)$ with $k=1, 2$,
as in \cite{Lon4}, the $\;\dm$-product of $M_1$ and $M_2$ is defined by
the following $2(m_1+m_2)\times 2(m_1+m_2)$ matrix $M_1\dm M_2$:
$$ M_1\dm M_2=\left(\matrix{A_1&  0&B_1&  0\cr
                               0&A_2&  0&B_2\cr
                             C_1&  0&D_1&  0\cr
                               0&C_2&  0&D_2\cr}\right). \nn$$  %\dm=\diamond
Denote by $M^{\dm k}$ the $k$-fold $\dm$-product $M\dm\cdots\dm M$. Note
that the $\dm$-product of any two symplectic matrices is symplectic. For any two
paths $\ga_j\in\P_{\tau}(2n_j)$ with $j=0$ and $1$, let
$\ga_0\dm\ga_1(t)= \ga_0(t)\dm\ga_1(t)$ for all $t\in [0,\tau]$.

A special path $\xi_n\in\P_{\tau}(2n)$ is defined by
\be \xi_n(t) = \left(\matrix{2-\frac{t}{\tau} & 0 \cr
                                             0 &  (2-\frac{t}{\tau})^{-1}\cr}\right)^{\dm n}
         \qquad {\rm for}\;0\le t\le \tau.  \lb{6.1}\ee
{\bf Definition 6.1.} (cf. \cite{Lon2}, \cite{Lon4}) {\it For any $\om\in\U$ and
$M\in \Sp(2n)$, define
\be  \nu_{\om}(M)=\dim_{\C}\ker_{\C}(M - \om I_{2n}).  \lb{6.2}\ee
For any $\tau>0$ and $\ga\in \P_{\tau}(2n)$, define
\be  \nu_{\om}(\ga)= \nu_{\om}(\ga(\tau)).  \lb{6.3}\ee

If $\ga\in\P_{\tau,\om}^{\ast}(2n)$, define
\be i_{\om}(\ga) = [\Sp(2n)_{\om}^0: \ga\ast\xi_n],  \lb{6.4}\ee
where the right hand side of (\ref{6.4}) is the usual homotopy intersection
number, and the orientation of $\ga\ast\xi_n$ is its positive time direction under
homotopy with fixed end points.

If $\ga\in\P_{\tau,\om}^0(2n)$, we let $\mathcal{F}(\ga)$
be the set of all open neighborhoods of $\ga$ in $\P_{\tau}(2n)$, and define
\be i_{\om}(\ga) = \sup_{U\in\mathcal{F}(\ga)}\inf\{i_{\om}(\beta)\,|\,
                       \beta\in U\cap\P_{\tau,\om}^{\ast}(2n)\}.
               \lb{6.5}\ee
Then
$$ (i_{\om}(\ga), \nu_{\om}(\ga)) \in \Z\times \{0,1,\ldots,2n\}, $$
is called the index function of $\ga$ at $\om$. }

Note that when $\om=1$, this index theory was introduced by
C. Conley-E. Zehnder in \cite{CoZ1} for the non-degenerate case with $n\ge 2$,
Y. Long-E. Zehnder in \cite{LZe1} for the non-degenerate case with $n=1$,
and Y. Long in \cite{Lon1} and C. Viterbo in \cite{Vit2} independently for
the degenerate case. The case for general $\om\in\U$ was defined by Y. Long
in \cite{Lon2} in order to study the index iteration theory (cf. \cite{Lon4}
for more details and references).

For any symplectic path $\ga\in\P_{\tau}(2n)$ and $m\in\N$,  we
define its $m$-th iteration $\ga^m:[0,m\tau]\to\Sp(2n)$ by
\be \ga^m(t) = \ga(t-j\tau)\ga(\tau)^j, \qquad
  \forall j\tau\leq t\leq (j+1)\tau,\;j=0,1,\ldots,m-1.
     \lb{6.6}\ee
We still denote the extended path on $[0,+\infty)$ by $\ga$.

{\bf Definition 6.2.} (cf. \cite{Lon2}, \cite{Lon4}) {\it For any $\ga\in\P_{\tau}(2n)$,
we define
\be (i(\ga,m), \nu(\ga,m)) = (i_1(\ga^m), \nu_1(\ga^m)), \qquad \forall m\in\N.
   \lb{6.7}\ee
The mean index $\hat{i}(\ga,m)$ per $m\tau$ for $m\in\N$ is defined by
\be \hat{i}(\ga,m) = \lim_{k\to +\infty}\frac{i(\ga,mk)}{k}. \lb{6.8}\ee
For any $M\in\Sp(2n)$ and $\om\in\U$, the {\it splitting numbers} $S_M^{\pm}(\om)$
of $M$ at $\om$ are defined by
\be S_M^{\pm}(\om)
     = \lim_{\ep\to 0^+}i_{\om\exp(\pm\sqrt{-1}\ep)}(\ga) - i_{\om}(\ga),
   \lb{6.9}\ee
for any path $\ga\in\P_{\tau}(2n)$ satisfying $\ga(\tau)=M$.}

For a given path $\gamma\in {\cal P}_{\tau}(2n)$ we consider to deform
it to a new path $\eta$ in ${\cal P}_{\tau}(2n)$ so that
\begin{equation}
i_1(\gamma^m)=i_1(\eta^m),\quad \nu_1(\gamma^m)=\nu_1(\eta^m), \quad
         \forall m\in {\bf N}, \label{6.10}
\end{equation}
and that $(i_1(\eta^m),\nu_1(\eta^m))$ is easy enough to compute. This
leads to finding homotopies $\delta:[0,1]\times[0,\tau]\to {\rm Sp}(2n)$
starting from $\gamma$ in ${\cal P}_{\tau}(2n)$ and keeping the end
points of the homotopy always stay in a certain suitably chosen maximal
subset of ${\rm Sp}(2n)$ so that (\ref{6.10}) always holds. In fact,  this
set was first introduced in \cite{Lon2} as the path connected component
$\Omega^0(M)$ containing $M=\gamma(\tau)$ of the set
\begin{eqnarray}
  \Omega(M)=\{N\in{\rm Sp}(2n)\,&|&\,\sigma(N)\cap{\bf U}=\sigma(M)\cap{\bf U}\;
{\rm and}\;  \nonumber\\
 &&\qquad \nu_{\lambda}(N)=\nu_{\lambda}(M)\;\forall\,
\lambda\in\sigma(M)\cap{\bf U}\}. \label{6.11}
\end{eqnarray}
Here $\Omega^0(M)$ is called the {\it homotopy component} of $M$ in
${\rm Sp}(2n)$.

In \cite{Lon2}-\cite{Lon4}, the following symplectic matrices were introduced
as { \it basic normal forms}:
\begin{eqnarray}
D(\lambda)=\left(\matrix{\lm & 0\cr
         0  & \lm^{-1}\cr}\right), &\quad& \lm=\pm 2,\lb{6.12}\\
N_1(\lm,b) = \left(\matrix{\lm & b\cr
         0  & \lm\cr}\right), &\quad& \lm=\pm 1, b=\pm1, 0, \lb{6.13}\\
R(\th)=\left(\matrix{\cos\th & -\sin\th\cr
        \sin\th  & \cos\th\cr}\right), &\quad& \th\in (0,\pi)\cup(\pi,2\pi),
                     \lb{6.14}\\
N_2(\om,b)= \left(\matrix{R(\th) & b\cr
              0 & R(\th)\cr}\right), &\quad& \th\in (0,\pi)\cup(\pi,2\pi),
                     \lb{6.15}\end{eqnarray}
where $b=\left(\matrix{b_1 & b_2\cr
               b_3 & b_4\cr}\right)$ with  $b_i\in\R$ and  $b_2\not=b_3$.

Splitting numbers possess the following properties:

{\bf Lemma 6.3.} (cf. \cite{Lon2} and Lemma 9.1.5 of \cite{Lon4}) {\it Splitting
numbers $S_M^{\pm}(\om)$ are well defined, i.e., they are independent of the choice
of the path $\ga\in\P_\tau(2n)$ satisfying $\ga(\tau)=M$ appeared in (\ref{6.9}).
For $\om\in\U$ and $M\in\Sp(2n)$, splitting numbers $S_N^{\pm}(\om)$ are constant
for all $N\in\Om^0(M)$. }

{\bf Lemma 6.4.} (cf. \cite{Lon2}, Lemma 9.1.5 and List 9.1.12 of \cite{Lon4})
{\it For $M\in\Sp(2n)$ and $\om\in\U$, there hold
\begin{eqnarray}
S_M^{\pm}(\om) &=& 0, \qquad {\it if}\;\;\om\not\in\sg(M).  \lb{6.16}\\
S_{N_1(1,a)}^+(1) &=& \left\{\matrix{1, &\quad {\rm if}\;\; a\ge 0, \cr
0, &\quad {\rm if}\;\; a< 0. \cr}\right. \lb{6.17}\eea

For any $M_i\in\Sp(2n_i)$ with $i=0$ and $1$, there holds }
\be S^{\pm}_{M_0\dm M_1}(\om) = S^{\pm}_{M_0}(\om) + S^{\pm}_{M_1}(\om),
    \qquad \forall\;\om\in\U. \lb{6.18}\ee

We have the following

{\bf Theorem 6.5.} (cf. \cite{Lon3} and Theorem 1.8.10 of \cite{Lon4}) {\it For
any $M\in\Sp(2n)$, there is a path $f:[0,1]\to\Om^0(M)$ such that $f(0)=M$ and
\be f(1) = M_1(\om_1)\dm\cdots\dm M_k(\om_k),  \lb{6.19}\ee
where each $M_i(\om_i)$ is a basic normal form as in (\ref{6.12})-(\ref{6.15})
for $1\leq i\leq k$.}

\subsection{Multiplicity and stability of closed characteristics}

Let $\Sigma\in\H(2n)$. Fix a constant $\alpha$ satisfying
$1<\alpha<2$ and define the Hamiltonian function $H_\alpha:\R^{2n}\to [0,+\infty)$ by
\be  H_\alpha(x) = j(x)^{\alpha}, \qquad \forall x\in\R^{2n},  \lb{6.20}\ee
where $j$ is the gauge function of $\Sigma$, i.e., $j(x)=\lm$ if $x=\lm y$ for
some $\lm>0$ and $y\in\Sigma$ when $x\in\R^{2n}\bs\{0\}$, and $j(0)=0$.

Then $H_{\alpha}\in C^1(\R^{2n},\R)\cap C^3(\R^{2n}\bs\{0\},\R)$ is strictly
convex and $\Sigma=H_{\alpha}^{-1}(1)$. It is well-known that the problem
(\ref{1.1}) is equivalent to the following given energy problem of the
Hamiltonian system
\bea
\cases{\dot{y}(t) &$= JH_{\alpha}^{\prime}(y(t)),
\quad H_{\alpha}(y(t)) =1, \qquad \forall t\in\R,$ \cr
y(\tau) &$= y(0).$ \cr}
\lb{6.21}\eea
Denote by $\T(\Sigma,\alpha)$ the set of all geometrically distinct
solutions $(\tau,y)$ of the problem (\ref{6.21}).  Note that elements
in $\T(\Sigma)$ defined in Section 1 and $\T(\Sigma,\alpha)$ are one to one
correspondent to each other.

Let $(\tau,y)\in \T(\Sigma,\alpha)$. The fundamental solution
$\gamma_y:[0,\tau]\to\Sp(2n)$ with $\gamma_y(0)=I_{2n}$ of the linearized
Hamiltonian system
\be \dot{\xi}(t) = JH_{\alpha}''(y(t))\xi(t), \qquad
       \forall t\in\R,  \lb{6.22}\ee
is called the {\it associated symplectic path} of $(\tau,y)$.
The eigenvalues of $\gamma_y(\tau)$ are called {\it Floquet
multipliers} of $(\tau,y)$.

For any $(\tau,y)\in\T(\Sigma,\alpha)$ and $m\in\N$, we define
its $m$-th iteration $y^m:\R/(m\tau\Z)\to\R^{2n}$ by
\be y^m(t) = y(t-j\tau), \qquad \forall j\tau\leq t\leq (j+1)\tau,
       \quad j=0,1,2,\ldots, m-1. \lb{6.23}\ee
We still denote by $y$ its extension to $[0,+\infty)$.

We define via Definition 6.2 the following
\bea  S^+(y) &=& S_{\ga_y(\tau)}^+(1),  \lb{6.24}\\
  (i(y,m), \nu(y,m)) &=& (i(\ga_y,m), \nu(\ga_y,m)),  \lb{6.25}\\
   \hat{i}(y,m) &=& \hat{i}(\ga_y,m),  \lb{6.26}\eea
for all $m\in\N$, where $\ga_y$ is the associated symplectic path of $(\tau,y)$.

We have the following result:

{\bf Theorem 6.6.} (cf. Theorem 15.1.1 of \cite{Lon4} and references there in)
{\it Suppose $(\tau,y)\in \T(\Sigma)$. Then we have
\be i(y^m)\equiv i(m\tau ,y)=i(y, m)-n,\quad \nu(y^m)\equiv\nu(m\tau, y)=\nu(y, m),
       \qquad \forall m\in\N, \lb{6.27}\ee
where $i(m\tau ,y)$ and $\nu(m\tau ,y)$ are given by Definition 3.14.
In particular, we have
\be \hat{i}(\tau, y)=\hat{i}(y,1),\lb{6.28}\ee
where $\hat{i}(\tau ,y)$ is given by Definition 3.14. Hence we denote it simply
by $\hat{i}(y)$. }

Now we can prove Theorem 1.1 as follows:

For $n\ge 2$ and $\Sg\in\H(2n)$ with $\,^{\#}\T(\Sg)<+\infty$,
using the index iteration theory developed by Y. Long and his
coworkers, specially the common index jump theorem (Theorem 4.3 of
\cite{LoZ1}, Theorem 11.2.1 of \cite{Lon4}), we obtain the following
estimate on the number $\,^{\#}\T(\Sg)$ by Theorem 5.1 of \cite{LoZ1} (Theorem
15.4.3 of \cite{Lon4}):
\be \,^{\#}\T(\Sg) \ge \varrho_n(\Sg). \lb{6.29}\ee
Here the invariant $\varrho_n(\Sg)$ is defined to be the minimum value of
$[\frac{i(x,1) + 2S^+(x) - \nu(x,1)+n}{2}]$ for all infinitely variationally
visible closed characteristic $(\tau,x)\in \T(\Sg)$ (cf. Definition
1.1 of \cite{LoZ1}, Definition 15.4.1 of \cite{Lon4}).
Specially we obtain
\be \varrho_n(\Sg) \ge
   \min\left\{\left[\frac{i(x,1) + 2S^+(x) - \nu(x,1)+n}{2}\right]\,
        \left|\frac{}{}\right.\,(\tau,x)\in\T(\Sg)\right\}.
   \lb{6.30}\ee

In the proof of Theorem 1.1 of \cite{LoZ1}, i.e., the following estimate
(15.5.21) on page 340 of \cite{Lon4} holds
$$ 2S^+(x) - \nu(x,1) \ge 1-p_+ \ge -1, $$
where $p_+$ counts the number of basic normal form $N_1(1,-1)$ appears in
the basic normal form decomposition of $\ga_x(\tau)$ in $\Om^0(\ga_x(\tau))$.
This estimate indicates that the worst case for getting a better estimate on
$\varrho_n(\Sg)$ happens when $p_+=2$ holds. Here we have used Theorem 6.5 and
Lemma 6.3.

Now if there are only two geometrically distinct closed characteristics on
$\Sg\subset \R^6$, together with (\ref{1.4}) and (\ref{6.29}), our
above Theorem 1.2 can be used to either kill at least one of the possible
$N_1(1,-1)$s or to derive a contradiction when there are two $N_1(1,-1)$s in
the decomposition of $\ga_j(\tau_j)$ in $\Om^0(\ga_j(\tau_j))$. Thus the conjecture
(\ref{1.2}) holds for $n=3$. More precisely we have

{\bf Proof of Theorem 1.1.} Assume the contrary, i.e., by \cite{EkH1} or
\cite{LoZ1} we assume $\,^{\#}\T(\Sg)= 2$ for some $\Sg\in\H(6)$. We use the
techniques in the index iteration theory developed by Y. Long and
his coauthors (cf. \cite{Lon4}), specially those techniques in the proof of
Theorem 5.1 of \cite{LoZ1} (cf. p.340 of \cite{Lon4}) to reach a contradiction.

Denote the two prime closed characteristics on $\Sg$ by $(\tau_j,\,y_j)$ with
the corresponding associated symplectic paths
$\ga_j\equiv \gamma_{y_j}:[0,\tau_j]\to\Sp(6)$ for $j=1,\,2$.
Then by Lemma 1.3 of \cite{LoZ1} or Lemma 15.2.4 of \cite{Lon4}, there exist
$P_j\in \Sp(6)$ and $M_j\in \Sp(4)$ such that
$\gamma_j(\tau_j)=P_j^{-1}(N_1(1,\,1)\diamond M_j)P_j$.
By our Theorem 1.2, we obtain the following identity:
\be \frac{\hat{\chi}(y_1)}{\hat{i}(y_1)}
         + \frac{\hat{\chi}(y_2)}{\hat{i}(y_2)}=\frac{1}{2}.  \lb{6.31}\ee
It is well known that $\hat{i}(y_j)>2$ for $j=1$ and $2$
(cf. Theorem 1.7.7 of \cite{Eke3} or Lemma 15.3.2 of \cite{Lon4}).

By Theorems 1.1 and 1.3 of \cite{LoZ1} (cf. Theorems 15.4.3 and 15.5.2 of
\cite{Lon4}) with $n=3$, we may assume that $y_1$ has irrational mean index
$\hat{i}(y_1)$. Next we continue our study in two cases.

{\bf Case 1}. {\it The average Euler characteristic $\hat{\chi}(y_1)\not=0$.}

In this case, by (\ref{6.31}) both $y_1$ and $y_2$ must possess irrational mean
indices. Hence by Theorem 8.3.1 and Corollary 8.3.2 of \cite{Lon4}, each $M_j$
can be connected to $R(\theta_j)\dm N_j$ within $\Om^0(M_j)$ for some
$\frac{\theta_j}{\pi}\notin\Q$ and $N_j\in \Sp(2)$. Now by Lemma 6.4 we have
\bea
&& S^+_{N_1(1,\,1)}(1)=\nu_1(N_1(1,\,1))=1,
    \quad S^+_{R(\theta_j)}(1)=\nu_1(R(\theta_j))=0, \lb{6.32}\\
&& 2S^+_{N_j}(1) - \nu_1(N_j) \ge -1. \lb{6.33}\eea
Thus we obtain by (\ref{6.10}), Lemma 6.3 and Lemma 6.4
\bea
&& 2S^+(y_j)-\nu(y_j,\,1) \nn\\
&&\qquad = 2S^+_{N_1(1,\,1)}(1)-\nu_1(N_1(1,\,1))
   +2S^+_{R(\theta_j)}(1)-\nu_1(R(\theta_j))+2S^+_{N_j}(1)-\nu_1(N_j) \nn\\
&&\qquad = 1 + 2S^+_{N_j}(1) - \nu_1(N_j)\nn\\
&&\qquad \ge 0, \lb{6.34}\eea
By Corollary 15.1.4 of \cite{Lon4} we have $i(y_j,\, 1)\ge 3$ for $j=1, 2$.
Therefore for $j=1, 2$, we obtain
\be i(y_j,\,1) + 2S^+(y_j) - \nu(y_j,\,1) \ge 3. \lb{6.35}\ee
Now by the estimates (\ref{6.29}) and (\ref{6.30}), we get a contradiction
\be  2= \;^{\#}\T(\Sg)\ge \varrho_3(\Sg) \ge 3, \lb{6.36}\ee
which completes the proof of Theorem 1.1 in Case 1.

{\bf Case 2}. {\it The average Euler characteristic $\hat{\chi}(y_1) =0$.}

In this case (\ref{6.31}) becomes
\be \frac{\hat{\chi}(y_2)}{\hat{i}(y_2)} = \frac{1}{2}.  \lb{6.37}\ee
Our above discussions in Case 1 can be applied to get (\ref{6.35}) for $j=1$.
Thus by Corollary 1.1 of \cite{LoZ1} (Theorem 15.4.4 of \cite{Lon4}), we
should get (\ref{6.36}) whenever (\ref{6.35}) holds for $j=2$. This yields a
contradiction.

Therefore now we assume that (\ref{6.35}) does not hold for $j=2$. Then as in
Case 1, we denote the basic normal form decomposition of $\ga_{y_2}(\tau_2)$ in
$\Om^0(\ga_{y_2}(\tau_2))$ by $N_1(1,1)\dm M_2$. By theorem 6.5, the $4\times 4$
matrix $M_2$ is either the $\diamond$-product of two matrices in
(\ref{6.12})-(\ref{6.14}) or one matrix in (\ref{6.15}). Therefore by Lemmas
6.3-6.4 and the first part of (\ref{6.32}), we have
$$ 2S^+(y_2) - \nu(y_2, 1)= 1 + S^+_{M_2}(1)-\nu_1(M_2)\ge 1-p_+ \ge -1, $$
where $p_+$ counts the number of basic normal form $N_1(1,-1)$ appears in $M_2$,
and both of the last two equalities hold simultaneously if and only if $p_+=2$,
i.e.,
\be M_2=N_1(1,-1)^{\dm 2}.  \lb{6.38}\ee
Because $i(y_2,1)\ge 3$, we obtain that the only case for which
(\ref{6.35}) and consequently Theorem 1.1 does not hold is when (\ref{6.38}) and
\be  i(y_2,1)=3  \lb{6.39}\ee
hold. Hence in the following, it suffices to derive a contradiction in this case.

Now by Theorem 8.3.1 of \cite{Lon4}, we obtain
\be i(y_2,m)=m(i(y_2,1)+1)-1=4m-1, \quad \nu(y_2,m) = 3, \quad \forall\,m\in\N.
        \lb{6.40}\ee
By Theorem 6.6,  we have
\be i(y_2)=i(y_2,1)-3=0, \quad i(y_2^2)=i(y_2,2)-3=4, \quad \hat{i}(y_2)=4.
   \lb{6.41}\ee
By Proposition 3.13, we obtain $K(y_2)=2$. By Remark 3.17 and (\ref{1.4})
with $n=3$ we obtain
\bea  \hat{\chi}(y_2)
&=& \frac{1}{K(y_2)}\sum_{1\le m\le 2\atop 0\le l\le 2}(-1)^{i(y_2^{m})+l}k_l(y_2^{m})
     \nn\\
&=& \frac{1}{2}(k_0(y_2) - k_1(y_2) + k_2(y_2) + k_0(y_2^2) - k_1(y_2^2) + k_2(y_2^2)) \nn\\
&\le& 1.  \lb{6.42}\eea
Now (\ref{6.37}), (\ref{6.41}) and (\ref{6.42}) yield a contradiction:
$$  \frac{1}{2} = \frac{\hat{\chi}(y_2)}{\hat{i}(y_2)} \le \frac{1}{4}, $$
which proves Theorem 1.1 in Case 2.

The proof of Theorem 1.1 is complete. \hfill\hb

{\bf Proof of Theorem 1.4.} Using notations in the above proof of Theorem 1.1,
we obtain (\ref{6.31}) for the two prime closed characteristics $(\tau_1,y_1)$
and $(\tau_2,y_2)$ on $\Sg$ with $M_j\in\Sp(2)$ and $M_1=R(\th_1)$ for some
$\th_1\in \R\bs \pi\Q$. Then $y_1$ is non-degenerate and then we obtain
$\hat{\chi}(y_1)\not= 0$ by (\ref{1.4}) and (\ref{3.58}). Therefore $\hat{i}(y_2)$
has to be irrational by (\ref{6.31}). \hfill\hb

{\bf Remark 6.7.} Using notations in the proof of Theorem 1.1, by Theorem
1.4 for $j=1$ and $2$ there exists $P_j\in\Sp(4)$ such that
$\ga_j(\tau_j)=P_j^{-1}(N_1(1,\,1)\dm R(\th_j))P_j$ for some
$\th_j\in (0,\pi)\cup(\pi,2\pi)$ with $\th_j/\pi\not\in \Q$. Then we obtain
$\nu(y_j^m)\equiv1$ and $i(y_j^m)\in 2\Z$,
specially $y_j^m$ are all non-degenerate for all $m\in\N$ and then
$K(y_j)=2$ for $j=1,2$. Thus we have $\hat{\chi}(y_j)=1$ for $j=1,2$
by (\ref{1.4}). Then together with (\ref{6.31}) we obtain that
$\hat{i}(y_1)/\hat{i}(y_2)$ is irrational. Therefore in this sense, $\Sg$
behaves like a weakly non-resonant ellipsoid (cf. \cite{Eke3}).

\medskip

\noindent {\bf Acknowledgements.} The authors sincerely thank Professor
Jianzhong Pan for helpful discussions with him on some topological
result related to Lemma 3.9 of this paper, and Professor Ivar Ekeland
for explaining to them in 2006 some results of the paper \cite{Eke1}.
Y. Long sincerely thanks Professor Claude Viterbo for explaining to him
certain details of the paper \cite{Vit1} in 2001. The authors sincerely
thank the referees for their careful readings of this paper and valuable
comments and suggestions.

\bibliographystyle{abbrv}

\bigskip

First submitted: April 23rd, 2006.

Revised version submitted: September 7th, 2006.

\end{document}